\documentclass[a4paper,10pt]{article}

\usepackage[dvips]{graphics,graphicx,epsfig}
\usepackage{fancybox}
\usepackage{latexsym}
\usepackage{amssymb}
\usepackage[all]{xy}
\usepackage{color}

\newcommand{\keywordsname}{Keywords}
\newenvironment{keywords}%
  {\small
    \list{}{\labelwidth0pt
      \leftmargin0pt \rightmargin\leftmargin
      \listparindent\parindent \itemindent0pt
      \parsep0pt
      }%
    \item[\hskip\labelsep\bfseries\keywordsname:]}{\endlist}

\def\Nset{\mathrm{I\!N}}
\def\Rset{\mathrm{I\!R}}

\def\rien{\rule{0pt}{0pt}}
\def\text#1{\mbox{#1}}
\def\counterfact{\Box\!\!\rightarrow}

\begin{document}


\title{D\emph{eterministic} m\emph{odal} B\emph{ayesian} L\emph{ogic}: derive the Bayesian within the modal logic $T$}

\author{\begin{tabular}{c}
{\bf Fr\'ed\'eric Dambreville}\\
D\'el\'egation G\'en\'erale pour l'Armement, DGA/CEP/GIP/SRO\\
16 Bis, Avenue Prieur de la C\^ote d'Or\\
F 94114, France\\[3pt]
{\tt http://email.FredericDambreville.com}\\
{\tt http://www.FredericDambreville.com}
\end{tabular}}

\date{}
\maketitle

\begin{abstract}
In this paper a conditional logic is defined and studied.
This conditional logic, DmBL, is constructed as close as possible to the Bayesian and is unrestricted, that is one is able to use any operator without restriction.
A notion of logical independence is also defined within the logic itself.
This logic is shown to be non trivial and is not reduced to classical propositions.
A model is constructed for the logic.
Completeness results are proved.
It is shown that any unconditioned probability can be extended to the whole logic DmBL.
The Bayesian is then recovered from the probabilistic DmBL.
At last, it is shown why DmBL is compliant with Lewis triviality.
\end{abstract}

\section{Introduction}
This work originates from a contribution to the book of Dezert and Smarandache \cite{dezert} about the DSmT, a theory related to the Dempster-Shafer Theory of Evidence.
Evidence Theory \cite{shafer} is a theoretical and practical tool for manipulating non deterministic informations.
In particular, it is used for fusing information.
Among other non deterministic theories, Evidence Theory is considered by some as a possible alternative to the Bayesian Probability.
\\[5pt]
A question then arises.
Why using a theory instead another theory?
This question implies another sub-question: what is the logic behind a non deterministic theory?
For example, when is known the logic behind the DSmT or the Bayesian Probability, then it may be possible to make a comparison between these two theories.
It has been shown in \cite{dezert,dambreville} that such a logic can be found for the DSmT (it was a multi-Modal logic).
This paper is the development of a similar work about the Bayesian rule, which beginnings are found in~\cite{dambreville}.
A future work will be certainly a joint logical study of these theories.
\\[5pt]
There are several existing logical characterizations of the Bayesian Probability.
First should be mentioned Cox's axiomatic derivation of the Bayesian probability (\cite{cox} but also \cite{debrucq1,debrucq2}).
On the other hand, the theory of Probabilistic Logic \cite{adams,pearl,nilss,paass,andersen} explains how probabilities and the probabilistic Bayesian rule could be an approach for reasoning with uncertain propositions.
An orthogonal viewpoint consists in defining a logic or an algebra of conditionals, and then consider their possible probabilistic extensions.
This is the domain of conditional logics\cite{Lewis:counterfactuals,giordano1, giordano2} or of Conditional Event Algebras \cite{goodman,goodman2, dubois,calabrese};
in particular, the CEAs described in \cite{goodman,goodman2} make possible some enrichment of the Bayesian rules. 
By the way, it is known that the construction of a conditional logic/algebra is strongly constrained by Lewis result \cite{lewis}, which has shown some critical issues related to the notion of conditional probability; refer also to \cite{hajek1,hajek2,fraassen}.
In particular, Lewis result implies strong hypotheses about the nature of the conditionals.
In most cases, the conditionals have to be constructed outside the space of unconditioned propositions.
Such constraints are apparent in the model constructed in \cite{goodman2}\,.
Whatever, an escape to Lewis, may be to restrict the construction of the conditionals to unconditioned propositions only~\cite{goodman2}\,.
Closure of the conditionals are proposed for several CEAs, but some undesired properties come out~\cite{goodman}.
\\[3pt]
In this paper, a conditional logic is defined, with the purpose to be closely related to a Bayesian reasoning (notice that \cite{goodman2} gives an example of conditional algebra very close to the Bayesian).
Moreover this logic is constructed so as to be closed for the conditionals.
This closure makes possible the definition and use of the notion of logical independence within the logic.
It is shown that the logical independence have several logical consequences; independence properties are derived.
At last, the Bayesian probability is recovered from these logical properties and the probabilistic extension of DmBL.
As a final result, a theorem is proved that guarantees the existence of an extension of any probability over classical proposition onto DmBL (Lewis result is thus avoided).
\\[3pt]
As will be shown along the paper, DmBL does not add notable properties to the already existing CEA\,, at the unconditioned level.
But its main interest is to illustrate the logical behavior of the conditional closure.
\\\\
Section~\ref{Section:Def:DmBL} is dedicated to the definition of the Deterministic modal Bayesian Logic.
The languages, axioms and rules are introduced, followed by an informal discussion of the definition.
\\[5pt]
In section~\ref{Section:Theorem:DmBL}, several theorems of the logic are derived.
A discussion about Lewis' negative result is opened here.
A purely logical viewpoint is considered in this section.
\\[5pt]
A model for DmBL is constructed in section~\ref{Section:Model:DmBL}.
A partial completeness theorem is derived.
\\[5pt]
The extension of probabilities over DmBL is investigated in section~\ref{Section:Proba:DmBL}.
The (probabilistic) Bayesian rule is recovered from the probabilistic extension of the Deterministic modal Bayesian Logic.
It is explained why our construction avoids Lewis' negative result.
\\[5pt]
The paper is then concluded.
\\[5pt]
Some comparisons of our model and logic with existing works are done after the definition, the logical theorems, the model construction (in appendix) and the probability extension.
\section{Definition of the logics}
\label{Section:Def:DmBL}
DmBL was first defined as DBL (Deterministic Bayesian Logic) without modality in a previous version of this document\cite{Dambreville:DBL}.
The non modal definition is really harsh and uneasy to handle/understand.
For this new modal definition, I have been inspired by the seminal work of Lewis \cite{Lewis:counterfactuals}, and also by more recent works of Laura Giordano and al, which use the modality for specifying the conditional behavior~\cite{giordano1,giordano2}.
In the following, the definition of DmBL is also based on an interaction with the modality.
As will be discussed later, the modality is mainly instrumental and the interaction is implicit in DmBL.
\\[5pt]
The D\emph{eterministic} m\emph{odal} B\emph{ayesian} L\emph{ogic} is defined subsequently.
It is an extension of the modal logic $T$, alias \emph{System} $T$, which definition is recalled.
In the definition of DmBL, the modal operator $\Box$ plays the role of a logical tool, for discriminating a proposition from the knowledge we have of it.
\paragraph{Language.}
Let $\Theta=\{\theta_i/i\in I\}$ be a set of atomic propositions.
\\[5pt]
The language $\mathcal{L}_C$ of the classical logic related to $\Theta$ is the smallest set such that:
\begin{itemize}
\item $\Theta\subset\mathcal{L}_C$ and $\bot,\top\in\mathcal{L}_C$\,,
\item $\neg\phi\in\mathcal{L}_C$\,, $\phi\rightarrow\psi\in\mathcal{L}_C$\,, $\phi\wedge\psi\in\mathcal{L}_C$ and $\phi\vee\psi\in\mathcal{L}_C$ for any $\phi,\psi\in\mathcal{L}_C$\,.
\end{itemize}
The language $\mathcal{L}_T$ of the modal logic $T$ related to $\Theta$ is the smallest set such that:
\begin{itemize}
\item $\Theta\subset\mathcal{L}_T$ and $\bot,\top\in\mathcal{L}_T$\,,
\item $\neg\phi\in\mathcal{L}_T$\,, $\Box\phi\in\mathcal{L}_T$ and $\Diamond\phi\in\mathcal{L}_T$ for any $\phi\in\mathcal{L}_T$\,,
\item $\phi\rightarrow\psi\in\mathcal{L}_T$\,, $\phi\wedge\psi\in\mathcal{L}_T$ and $\phi\vee\psi\in\mathcal{L}_T$ for any $\phi,\psi\in\mathcal{L}_T$\,.
\end{itemize}
The language $\mathcal{L}$ of the D\emph{eterministic} m\emph{odal} B\emph{ayesian} L\emph{ogic} related to $\Theta$ is the smallest set such that:
\begin{itemize}
\item $\Theta\subset\mathcal{L}$ and $\bot,\top\in\mathcal{L}$\,,
\item $\neg\phi\in\mathcal{L}$\,, $\Box\phi\in\mathcal{L}$ and $\Diamond\phi\in\mathcal{L}$ for any $\phi\in\mathcal{L}$\,,
\item $\phi\rightarrow\psi\in\mathcal{L}$\,, $\phi\wedge\psi\in\mathcal{L}$\,, $\phi\vee\psi\in\mathcal{L}$\,, $(\psi|\phi)\in\mathcal{L}$  and $\phi\times\psi\in\mathcal{L}$ for any $\phi,\psi\in\mathcal{L}$\,.
\end{itemize}
The operator $\times$ describes the \emph{logical} independence between propositions.
The operators $\times$ and $(|)$ are conjointly defined.
In the sequel, the notation $\phi\leftrightarrow\psi$ is an abreviation for $(\phi\rightarrow\psi)\wedge(\psi\rightarrow\phi)$\,.
Moreover $\phi\equiv\psi$ means $\vdash\phi\leftrightarrow\psi$\,.
\paragraph{Rules and axioms.}
The \emph{classical Logic $C$} is characterized by the \emph{Modus ponens} and the classical axioms $c\ast$ described subsequently.
\\
The \emph{modal Logic $T$} is characterized by the \emph{Modus ponens}, the classical axioms $c\ast$ and the modal rule/axioms $m\ast$ described subsequently (\emph{c.f.} also \cite{blackburn}).
\\
The D\emph{eterministic} m\emph{odal} B\emph{ayesian} L\emph{ogic}\,, \emph{i.e.} DmBL, is characterized by the \emph{Modus ponens}, the classical axioms $c\ast$, the modal rule/axioms $m\ast$ and the Bayesian axioms $b\ast$\,:
\begin{description}
\item[\rien$\quad$c1.] $\vdash\top$\,,
\item[\rien$\quad$c2.] $\vdash \phi\rightarrow(\psi\rightarrow\phi)$\,,
\item[\rien$\quad$c3.] $\vdash (\eta\rightarrow(\phi\rightarrow\psi))\rightarrow((\eta\rightarrow\phi)\rightarrow(\eta\rightarrow\psi))$\,,
\item[\rien$\quad$c4.] $\vdash (\neg\phi\rightarrow\neg\psi)\rightarrow((\neg\phi\rightarrow\psi)\rightarrow\phi)$\,,
\item[\rien$\quad$c5.] $\bot\equiv\neg\top$\,,
\item[\rien$\quad$c6.] $\phi\rightarrow\psi\equiv\neg\phi\vee\psi$\,,
\item[\rien$\quad$c7.] $\phi\wedge\psi\equiv\neg(\neg\phi\vee\neg\psi)$\,,
\item[\rien$\quad$Modus ponens.]$\vdash\phi$ and $\vdash\phi\rightarrow\psi$ implies $\vdash\psi$\,,
\item[\rien$\quad$m1.] $\vdash\phi$ implies $\vdash \Box\phi$\,,
\item[\rien$\quad$m2.] $\vdash\Box(\phi\rightarrow\psi)\rightarrow(\Box\phi\rightarrow\Box\psi)$\,,
\item[\rien$\quad$m3.] $\vdash\Box\phi\rightarrow\phi$\,,
\item[\rien$\quad$m4.] $\Diamond\phi\equiv\neg\Box\neg\phi$\,,
\item[\rien$\quad$b1.] $\vdash\Box(\phi\rightarrow\psi)\rightarrow \bigl(\Box\neg\phi \vee \Box(\psi|\phi)\bigr)$\,,
\item[\rien$\quad$b2.] $\vdash(\psi\rightarrow\eta|\phi)\rightarrow\bigl((\psi|\phi)\rightarrow(\eta|\phi)\bigr)$\,,
\item[\rien$\quad$b3.] $\vdash(\psi|\phi)\rightarrow(\phi\rightarrow\psi)$\,,
\item[\rien$\quad$b4.] $\neg(\neg\psi|\phi)\equiv(\psi|\phi)$\,,
\item[\rien$\quad$b5 (definition of independence).] $\psi\times\phi\equiv\Box\bigl((\psi|\phi)\leftrightarrow\psi\bigr)$\,,
\item[\rien$\quad$b6.] \emph{($\times$ is symmetric)}~:
$\psi\times\phi\equiv\phi\times\psi$\,,
\item[\rien$\quad$Customized axioms.] Axioms introduced for the specification of a particular sub-theory of DmBL.
In particular the customized axioms may be used to characterize the properties of the atomic propositions: axioms like $\vdash\theta_i$\,,$\;$\ $\vdash\neg\Box\theta_i$ and $\vdash\theta_i\times\theta_j$ are possible.
\end{description}
Customized axioms will not be considered in this paper.
\\[5pt]
DmBL$_\ast$, a weakened version of DmBL, is defined by replacing $b6$ by the alternative axioms:
$$\begin{array}{@{}l@{}}
\mbox{\bf b6.weak.A. } \psi\times\neg\phi\equiv\psi\times\phi\;,
\\[3pt]
\mbox{\bf b6.weak.B. } \vdash\Box(\psi\leftrightarrow\eta)\rightarrow\Box\bigl((\phi|\psi)\leftrightarrow(\phi|\eta)\bigr)\;.
\end{array}$$
The symbols $\vdash_C$, $\vdash_T$ and $\vdash$ will be used to denote a proof in $C$, $T$ or DmBL/DmBL$_\ast$ respectively.
\paragraph{Discussion.}
First at all, let us discuss the role of the modality.
Notice that the modality is only instrumental for us: there exists a non modal definition of DmBL, called DBL, defined in a previous version~\cite{Dambreville:DBL}.
This property will be confirmed in section~\ref{Section:Model:DmBL}, when is constructed a Kripke model of DmBL$_\ast$ complete for the non modal propositions (including the conditionals).
This Kripke model is based on a full accessibility relation (\emph{i.e.} $R=W^2$); full accessibility indeed reduces the modality to something trivially contrasted\,.
\\[3pt]
Now, the role of $\Box$ here is to infer a modal-like definition of the conditional $(|)$\,.
More precisely, the fundamental axiom b1: $\vdash\Box(\phi\rightarrow\psi)\rightarrow \bigl(\Box\neg\phi \vee \Box(\psi|\phi)\bigr)$\,, should be compared to its non modal version in DBL~\cite{Dambreville:DBL}\,, that is 
$$
\mathbf{b1.old.} \vdash\phi\rightarrow\psi\quad\mbox{implies}\quad \vdash\neg\phi\mbox{ OR }\vdash(\psi|\phi)\;,
$$
which is clearly a modal-like definition.
But b1.old is really uneasy to handle (what is a proof with such multivalued rule?), while working at the modal level with b1 is really quiet.
Now, it is interesting to compare our axiomatization with previously existing axiomatization of the conditionals.
\\[5pt]
Typically, we will consider the axioms of VCU~\cite{Lewis:counterfactuals}\,, p.132.
I guess that this comparison will make apparent some specific differences of DmBL.
Notice that this comparison will also rely on logical theorems derived subsequently in section~\ref{Section:Theorem:DmBL}\,, or from rather easy deductions:\\
{\small
(Ax.1) $\phi\counterfact\phi$ has a partial counterpart in DmBL, \emph{i.e.} $\vdash\Box\neg\phi\vee\Box(\phi|\phi)$ (theorem).\\
(Ax.2) $(\neg\phi\counterfact\phi)\rightarrow(\psi\counterfact\phi)$ becomes $\vdash(\phi|\neg\phi)\rightarrow(\phi|\psi)$ (derived from theorem.)\\
(Ax.3) $(\phi\counterfact\neg\psi)\vee(((\phi\wedge\psi)\counterfact\xi)\leftrightarrow(\phi\counterfact(\psi\rightarrow\xi)))$ has no obvious counterpart in DmBL.\\
(Ax.4) $(\phi\counterfact\psi)\rightarrow(\phi\rightarrow\psi)$ is related to b3\\
(Ax.5) $(\phi\wedge\psi)\rightarrow(\phi\counterfact\psi)$ becomes $\vdash(\phi\wedge\psi)\rightarrow(\psi|\phi)$ which is a part of the inference theorem.\\
(Ax.6) $(\neg\phi\counterfact\phi)\rightarrow\bigl(\neg(\neg\phi\counterfact\phi)\counterfact(\neg\phi\counterfact\phi)\bigr)$
becomes
$\vdash(\phi|\neg\phi)\rightarrow\bigl((\phi|\neg\phi)\big|\neg(\phi|\neg\phi)\bigr)$
(derived from theorem.)\\
(CR) Counterfactual rule.
This rule is a multiple-task rule!
First, it allows the introduction of tautologies inside a conditional,
secondly, it implies some linearity of the conditional with $\wedge$\,:
\\\rien\hspace{50pt}Being proved $(\xi_1\wedge\dots\wedge\xi_n)\rightarrow\psi$\,, it is proved $((\phi\counterfact\xi_1)\wedge\dots\wedge(\phi\counterfact\xi_n))\rightarrow(\phi\counterfact\psi)$\,.\\
This rule becomes in DmBL:
\\\rien\hspace{50pt}$\vdash(\xi_1\wedge\dots\wedge\xi_n)\rightarrow\psi$ implies $\vdash((\xi_1|\phi)\wedge\dots\wedge(\xi_n|\phi))\rightarrow(\psi|\phi)$\,.
\\[5pt]
As a first conclusion, it happens that Ax.2, Ax.4, Ax.5, Ax.6 and CR are recovered in DmBL, Ax.1 is partially recovered, and Ax.3 is not derived from DmBL.
In fact, the axiom Ax.3
is replaced in DmBL by a strong axiom dedicated to the negation, \emph{i.e.} b4.
It will be seen that the weakening of Ax.1 is forced by the axiom b4 of DmBL.
Otherwise, some conflicts will focus on the empty proposition $\bot$.
\\[3pt]
Now, let discuss precisely about the specific axioms of DmBL.
\\\\
In VCU, there is a rule, CR, for the introduction of tautologies inside a conditional.
DmBL contains this rule.
But in DmBL, the introduction of ``tautologies'' (in fact $\Box$-ed propositions) in the conditional is done more specifically by b1.
\begin{itemize}
\item b1 manages the ``tautologies'' more precisely.
While CR (simplified) maps the tautology $\xi\rightarrow\psi$ onto the tautology  $(\phi\counterfact\xi)\rightarrow(\phi\counterfact\psi)$\,,
b1 maps $\Box(\phi\rightarrow\psi)$ onto $\Box(\psi|\phi)$\,: \emph{i.e.} the proposition $\psi$ could be considered within the range $\phi$,
\item But the mapping of b1 is \emph{always} weakened by the output $\Box\neg\phi$\,.
\vspace{3pt}
\end{itemize}
Axiom b3 is the same as Ax.4.\\[5pt]
Axiom b2 has no equivalent in VCU.
However, b2 could be deduced from b4 and $(\psi\wedge\eta|\phi)\equiv(\psi|\phi)\wedge(\eta|\phi)$\,.
Now, it is proved $(\phi\counterfact(\psi\wedge\eta))\leftrightarrow((\phi\counterfact\psi)\wedge(\phi\counterfact\eta))$ in VCU.
In other words, b2 does not make a true difference with VCU.
\underline{The crucial difference thus comes with b4.}\\[5pt]
Axioms like b4 does not exist in VCU.
This axiom changes fundamentally the nature of the conditional: the sub-universe associated to a conditional becomes classical.
By classical, we mean that $(\psi\wedge\eta|\phi)\equiv(\psi|\phi)\wedge(\eta|\phi)$ and 
$(\neg\psi|\phi)\equiv\neg(\psi|\phi)$
(these properties are proved in the next section).
\\[3pt]
Now, it is remembered that axioms like b4 have not been included in VCU, because it implied some conflict with the other axioms/rule.
This conflict is particularly focused over the empty proposition $\bot$\,.
More precisely, it is deduced $\bot\counterfact\bot$ by Ax.1.
By CR and $\bot\rightarrow\top$, it is deduced $(\bot\counterfact\bot)\rightarrow(\bot\counterfact\top)$\,, and finally $(\bot\counterfact\bot)\wedge(\bot\counterfact\top)$\,.
If we apply a rule like b4, it is obtained $(\bot\counterfact\bot)\equiv\neg(\bot\counterfact\top)$ which contradicts $(\bot\counterfact\bot)\wedge(\bot\counterfact\top)$\,. 
\\[3pt]
In order to avoid such conflict in DmBL, the $\bot$ is considered as an exception in the rule b1, which results in the weakening of the mapping b1 by the proposition $\Box\neg\phi$\,.
\\[5pt]
Axiom b5 and b6 introduce quite new notions, in comparison with VCU.
But I do not have elements for comparison.
It is possible that b5 and b6 do not imply fundamental differences with VCU.
}
\\\\
Conclusion: in comparison with VCU, DmBL improves the manipulation of the negation within the conditional, and thus implies a classical behavior of the sub-universes related to a conditioning (\emph{c.f.} \emph{Sub-universes are classical} in next section).
As a consequence, it is necessary to introduce a relaxation in comparison with the rules of VCU, in particular around the empty proposition: $\bot$ should be considered as a singularity.
By manipulating the negation that way, DmBL is closer to ``Bayesian'' algebras (like for example the Conditional Events Algebra defined in~\cite{goodman2}; we will discuss this model later) than could be the VCU\,.
Moreover, DmBL is not restricted in the use of the conditionals (contrary to the CEA), and a notion of logical independence is constructed within the logic.
The following section studies the logical consequences of such axiomatization.
\section{Logical theorems}
\label{Section:Theorem:DmBL}
Notice first that DmBL/DmBL$_\ast$ contains the classical and the T-system tautologies.
The properties of classical logic and of the T-system will be assumed without proof subsequently.
Since we are studying both DmBL and DmBL$_\ast$, it is indicated in bracket the possibly needed axioms $b6\ast$.
\paragraph{The full universe} $\vdash\Box\phi\rightarrow(\psi\times\phi)$\,.
In particular $(\psi|\top)\equiv\psi$\,.\\[5pt]
Interpretation: a \emph{full} proposition is independent with any other proposition and its sub-universe is the whole universe.
\begin{description}
\item[Proof.]
From axiom b3, it comes $\vdash(\psi|\phi)\rightarrow(\phi\rightarrow\psi)$ and  $\vdash(\neg\psi|\phi)\rightarrow(\phi\rightarrow\neg\psi)$\,.\\
Then $\vdash\phi\rightarrow\bigl((\psi|\phi)\rightarrow\psi\bigr)$ and  $\vdash\phi\rightarrow\bigl((\neg\psi|\phi)\rightarrow\neg\psi\bigr)$\,.\\
Applying b4 yields $\vdash\phi\rightarrow\bigl((\psi|\phi)\leftrightarrow\psi\bigr)$.\\
It follows $\vdash\Box\phi\rightarrow\Box\bigl((\psi|\phi)\leftrightarrow\psi\bigr)$\,.
\\[5pt]
The remaining proof is obvious.
\item[$\Box\Box\Box$]\rien
\end{description}
\paragraph{Axioms order.} Axiom b6 implies b6.weak.A.
\begin{description}
\item[Proof.]
From b6, it is deduced $\psi\times\neg\phi\equiv\neg\phi\times\psi$\,.\\
Now $\neg\phi\times\psi\equiv\Box\bigl((\neg\phi|\psi)\leftrightarrow\neg\phi\bigr)\equiv\Box\bigl(\neg(\phi|\psi)\leftrightarrow\neg\phi\bigr)\equiv
\Box\bigl((\phi|\psi)\leftrightarrow\phi\bigr)\equiv\phi\times\psi$ by b4.\\
By applying b6 again, it comes $\psi\times\neg\phi\equiv\psi\times\phi$\,.
\item[$\Box\Box\Box$]\rien
\end{description}
\paragraph{The empty universe [b6.weak.A]} $\vdash\Box\neg\phi\rightarrow(\psi\times\phi)$\,.
In particular $(\psi|\bot)\equiv\psi$\,.
\begin{description}
\item[Proof.]
$\vdash\Box\neg\phi\rightarrow(\psi\times\neg\phi)$
implies 
$\vdash\Box\neg\phi\rightarrow(\psi\times\phi)$ (b6.weak.A).\\[5pt]
The remaining proof is obvious.
\item[$\Box\Box\Box$]\rien
\end{description}
\paragraph{Left equivalences.}
$\vdash\Box(\psi\leftrightarrow\eta)\rightarrow\Bigl(\Box\neg\phi\vee\Box\bigl((\psi|\phi)\leftrightarrow(\eta|\phi)\bigr)\Bigr)$\,.
\begin{description}
\item[Proof.]
From $\vdash (\psi\rightarrow\eta)\rightarrow\bigl(\phi\rightarrow(\psi\rightarrow\eta)\bigr)$\,, it is deduced $\vdash\Box(\psi\rightarrow\eta)\rightarrow\Box\bigl(\phi\rightarrow(\psi\rightarrow\eta)\bigr)$\,.\\
By applying axiom b1, it comes $\vdash\Box(\psi\rightarrow\eta)\rightarrow\bigl(\Box\neg\phi\vee\Box(\psi\rightarrow\eta|\phi)\bigr)$\,.\\
Then axiom b2 implies $\vdash\Box(\psi\rightarrow\eta)\rightarrow\Bigl(\Box\neg\phi\vee\Box\bigl((\psi|\phi)\rightarrow(\eta|\phi)\bigr)\Bigr)$.\\
Since $\psi$ and $\eta$ are exchangeable, the theorem is deduced.
\item[$\Box\Box\Box$]\rien
\end{description}
\emph{Corollary [b6.weak.A].} $\vdash\Box(\psi\leftrightarrow\eta)\rightarrow\Box\bigl((\psi|\phi)\leftrightarrow(\eta|\phi)\bigr)$.
\begin{description}
\item[Proof.]
It has been proved $\vdash\Box\neg\phi\rightarrow(\psi\times\phi)$, or equivalently $\vdash\Box\neg\phi\rightarrow\Box\bigl((\psi|\phi)\leftrightarrow\psi\bigr)$.\\
Of course, also holds $\vdash\Box\neg\phi\rightarrow\Box\bigl((\eta|\phi)\leftrightarrow\eta\bigr)$.\\
Then $\vdash\bigl(\Box\neg\phi\wedge\Box(\psi\leftrightarrow\eta)\bigl)\rightarrow\Box\bigl((\psi|\phi)\leftrightarrow(\eta|\phi)\bigr)$.\\
The corollary is then deduced from the main proposition.
\item[$\Box\Box\Box$]\rien
\end{description}
\emph{Corollary 2 [b6.weak.A].} $\psi\equiv\eta$ implies $(\psi|\phi)\equiv(\eta|\phi)$.
\\
Proof is immediate from corollary.
\paragraph{Sub-universes are classical [b6.weak.A].}
\begin{itemize}
\item $(\neg\psi|\phi)\equiv\neg(\psi|\phi)$\,,
\item $(\psi\wedge\eta|\phi)\equiv(\psi|\phi)\wedge(\eta|\phi)$\,,
\item $(\psi\vee\eta|\phi)\equiv(\psi|\phi)\vee(\eta|\phi)$\,,
\item $(\psi\rightarrow\eta|\phi)\equiv(\psi|\phi)\rightarrow(\eta|\phi)$\,.
\end{itemize}
\begin{description}
\item[Proof.]
The first theorem is a consequence of axiom b4.
\vspace{5pt}\\
From axiom b2\,, it is deduced $\vdash(\psi\rightarrow\neg\eta|\phi)\rightarrow\bigl((\psi|\phi)\rightarrow(\neg\eta|\phi)\bigr)$\,.\\
It is deduced $\vdash\bigl((\psi|\phi)\wedge\neg(\neg\eta|\phi)\bigr)\rightarrow\neg(\psi\rightarrow\neg\eta|\phi)$\,.\\
Applying b4, it comes $\vdash\bigl((\psi|\phi)\wedge(\eta|\phi)\bigr)\rightarrow(\psi\wedge\eta|\phi)$.
\\[5pt]
Now $\vdash\phi\rightarrow\bigl((\psi\wedge\eta)\rightarrow\psi\bigr)$ and b1 imply $\vdash\Box\neg\phi\vee\Box\bigl((\psi\wedge\eta)\rightarrow\psi\big|\phi\bigr)$.\\
By b2 it is deduced $\vdash\Box\neg\phi\vee\Box\bigl((\psi\wedge\eta|\phi)\rightarrow(\psi|\phi)\bigr)$.\\
It is similarly proved $\vdash\Box\neg\phi\vee\Box\bigl((\psi\wedge\eta|\phi)\rightarrow(\eta|\phi)\bigr)$.\\
At last $\vdash\Box\neg\phi\vee\Box\Bigl((\psi\wedge\eta|\phi)\rightarrow\bigl((\psi|\phi)\wedge(\eta|\phi)\bigr)\Bigr)$.\\[3pt]
Now it has been shown $\vdash\Box\neg\phi\rightarrow\Box\bigl((\Xi|\phi)\leftrightarrow\Xi\bigr)$, and considering $\Xi=\psi$, $\eta$ or $\phi\wedge\eta$, it is implied $\vdash\Box\neg\phi\rightarrow\Box\Bigl((\psi\wedge\eta|\phi)\rightarrow\bigl((\psi|\phi)\wedge(\eta|\phi)\bigr)\Bigr)$.\\
At last $\vdash\Box\Bigl((\psi\wedge\eta|\phi)\rightarrow\bigl((\psi|\phi)\wedge(\eta|\phi)\bigr)\Bigr)$ and then $\vdash(\psi\wedge\eta|\phi)\rightarrow\bigl((\psi|\phi)\wedge(\eta|\phi)\bigr)$.\\
The second theorem is then proved.
~\vspace{5pt}\\
Third theorem is a consequence of the first and second theorems.
~\vspace{5pt}\\
Last theorem is a consequence of the first and third theorems.
\item[$\Box\Box\Box$]\rien
\end{description}
\paragraph{Evaluating $(\top|\cdot)$ and $(\bot|\cdot)$ [b6.weak.A].} Is proved $\vdash\Box\psi\rightarrow\Box(\psi|\phi)$\,.
In particular $(\top|\phi)\equiv\top$ and $(\bot|\phi)\equiv\bot$\,.
\begin{description}
\item[Proof.]
From $\vdash\Box\psi\rightarrow\Box(\phi\rightarrow\psi)$ and b1, it comes $\vdash\Box\psi\rightarrow\bigl(\Box\neg\phi\vee\Box(\psi|\phi)\bigr)$\,.\\
Now $\vdash\Box\neg\phi\rightarrow\Box\bigl((\psi|\phi)\leftrightarrow\psi\bigr)$\,,
and consequently $\vdash\Box\psi\rightarrow\Box(\psi|\phi)$\,.
%
\item[$\Box\Box\Box$]\rien
\end{description}
\paragraph{Inference property.}
$(\psi|\phi)\wedge\phi\equiv\phi\wedge\psi$\,.
\begin{description}
\item[Proof.]
From b3 it comes $\vdash(\neg\psi|\phi)\rightarrow(\phi\rightarrow\neg\psi)$\,.\\
Then $\vdash\neg(\phi\rightarrow\neg\psi)\rightarrow\neg(\neg\psi|\phi)$ and 
$\vdash(\phi\wedge\psi)\rightarrow(\psi|\phi)$\,.\\
At last $\vdash(\phi\wedge\psi)\rightarrow\bigl((\psi|\phi)\wedge\phi\bigr)$\,.\\[5pt]
Conversely $\vdash(\psi|\phi)\rightarrow(\phi\rightarrow\psi)$ implies $\vdash\bigl((\psi|\phi)\wedge\phi\bigr)\rightarrow\bigl((\phi\rightarrow\psi)\wedge\phi\bigr)$\,.\\
Since $(\phi\rightarrow\psi)\wedge\phi\equiv\phi\wedge\psi$\,, the converse is proved.
\item[$\Box\Box\Box$]
\end{description}
\paragraph{Introspection.} $\vdash\Box\neg\phi\vee\Box(\phi|\phi)$\,.\\[5pt]
Interpretation: a non empty proposition sees itself as ever true. \\
Notice that this property is compliant with $(\bot|\bot)\equiv\bot$\,.
\begin{description}
\item[Proof.]
Obvious from $\vdash \phi\rightarrow\phi$ and b1\,.
\item[$\Box\Box\Box$]\rien
\end{description}
\paragraph{Inter-independence [b6.weak.A].} $\vdash(\psi|\phi)\times\phi$\,.\\[5pt]
Interpretation: a proposition is independent of its sub-universe.
\begin{description}
\item[Proof.]It is proved:
$$
\bigl((\psi|\phi)\big|\phi\bigr)\wedge(\phi|\phi)
\equiv\bigl((\psi|\phi)\wedge\phi\big|\phi\bigr)\equiv(\phi\wedge\psi|\phi)\equiv(\psi|\phi)\wedge(\phi|\phi)\;.
$$
As a consequence $\vdash(\phi|\phi)\rightarrow\Bigl(\bigl((\psi|\phi)\big|\phi\bigr)\leftrightarrow(\psi|\phi)\Bigr)$\,.\\
Then $\vdash\Box(\phi|\phi)\rightarrow\bigl((\psi|\phi)\times\phi\bigr)$\,.\\
Now $\vdash\Box\neg\phi\vee\Box(\phi|\phi)$ and $\vdash\Box\neg\phi\rightarrow\bigl((\psi|\phi)\times\phi\bigr)$ from previous results.\\
At last $\vdash(\psi|\phi)\times\phi$
\item[$\Box\Box\Box$]\rien
\end{description}
\paragraph{Independence invariance [b6.weak.A].}
$$
\begin{array}{@{}l@{}}
\vdash (\psi\times\phi)\rightarrow(\neg\psi\times\phi)\;,
\\
\vdash \bigl((\psi\times\phi)\wedge(\eta\times\phi)\bigr)\rightarrow\bigl((\psi\wedge\eta)\times\phi\bigr)\;,
\\
\vdash \Box(\psi\leftrightarrow\eta)\rightarrow\bigl((\psi\times\phi)\leftrightarrow(\eta\times\phi)\bigr)\;.
\end{array}
$$
\begin{description}
\item[Proof.]
First theorem comes from the deduction:
$$
(\psi\times\phi)\equiv\Box\bigl((\psi|\phi)\leftrightarrow\psi\bigr)\equiv\Box\bigl(\neg(\psi|\phi)\leftrightarrow\neg\psi\bigr)\equiv\Box\bigl((\neg\psi|\phi)\leftrightarrow\neg\psi\bigr)\equiv(\neg\psi\times\phi)\;.
$$
The second theorem is also derived from similar deductions:
$$
\vdash\Box\Bigl(\bigl((\psi|\phi)\leftrightarrow\psi\bigr)\wedge\bigl((\eta|\phi)\leftrightarrow\eta\bigr)\Bigr)
\rightarrow
\Box\Bigl(\bigl((\psi|\phi)\wedge(\eta|\phi)\bigr)\leftrightarrow(\psi\wedge\eta)\Bigr)\vspace{-5pt}
$$
and then $\vdash\Box\Bigl(\bigl((\psi|\phi)\leftrightarrow\psi\bigr)\wedge\bigl((\eta|\phi)\leftrightarrow\eta\bigr)\Bigr)
\rightarrow
\Box\bigl((\psi\wedge\eta|\phi)\bigr)\leftrightarrow(\psi\wedge\eta)\bigr)$\,.\vspace{5pt}
\\
Now, let prove the third theorem.\\
The \emph{Left equivalence} property implies $\vdash\Box(\psi\leftrightarrow\eta)\rightarrow\Box\Bigl(\bigl((\psi|\phi)\leftrightarrow(\eta|\phi)\bigr)\wedge(\psi\leftrightarrow\eta)\Bigr)$.\\
Since $\vdash\bigl((\alpha\leftrightarrow\beta)\wedge(\gamma\leftrightarrow\delta)\bigr)\rightarrow\bigl((\alpha\leftrightarrow\gamma)\leftrightarrow(\beta\leftrightarrow\delta)\bigr)$, it is deduced\\ \rien\hspace{150pt}$\vdash\Box(\psi\leftrightarrow\eta)\rightarrow\Box\Bigl(\bigl((\psi|\phi)\leftrightarrow\psi\bigr)\leftrightarrow\bigl((\eta|\phi)\leftrightarrow\eta\bigr)\Bigr)$.\\
At last $\vdash\Box(\psi\leftrightarrow\eta)\rightarrow\Bigl(\Box\bigl((\psi|\phi)\leftrightarrow\psi\bigr)\leftrightarrow\Box\bigl((\eta|\phi)\leftrightarrow\eta\bigr)\Bigr)$.
\item[$\Box\Box\Box$]\rien
\end{description}
\paragraph{Narcissistic independence.} $\vdash(\phi\times\phi)\rightarrow(\Box\neg\phi\vee\Box\phi)$\,.
\\[5pt]
Interpretation: a propositions independent with itself is either full or empty.
\begin{description}
\item[Proof.]\rien\\
From $\vdash\phi\rightarrow\phi$ it is deduced $\vdash\Box\neg\phi\vee\Box(\phi|\phi)$\,.\\
From b5 it is derived $\vdash(\phi\times\phi)\rightarrow\Box\bigl((\phi|\phi)\rightarrow\phi\bigr)$ and then $\vdash(\phi\times\phi)\rightarrow\bigl(\Box(\phi|\phi)\rightarrow\Box\phi\bigr)$\,.\\
It is thus deduced $\vdash(\phi\times\phi)\rightarrow(\Box\neg\phi\vee\Box\phi)$\,.
\item[$\Box\Box\Box$]\rien
\end{description}
\paragraph{Independence and proof [b6.weak.A].}
$\vdash(\psi\times\phi)\rightarrow\bigl(\Box(\phi\vee\psi)\rightarrow(\Box\phi\vee\Box\psi)\bigr)$\,.\\[5pt]
Interpretation:
when propositions are independent and their disjunctions are sure, then at least one proposition is sure.
\begin{description}
\item[Proof.]
Combining $\vdash(\phi\vee\psi)\rightarrow(\neg\phi\rightarrow\psi)$ with b1 implies $\vdash\Box(\phi\vee\psi)\rightarrow\bigl(\Box\phi\vee\Box(\psi|\neg\phi)\bigr)$\,.\\
From b6.weak.A, it comes $\vdash(\psi\times\phi)\rightarrow\Box\bigl((\psi|\neg\phi)\leftrightarrow\psi\bigr)$\,.\\
As a consequence $\vdash(\phi\times\psi)\rightarrow\bigl(\Box(\phi\vee\psi)\rightarrow(\Box\phi\vee\Box\psi)\bigr)$\,.
\item[$\Box\Box\Box$]\rien
\end{description}
\paragraph{Independence and regularity [b6.weak.A].}
$$
\vdash\bigl((\phi\times\eta)\wedge(\psi\times\eta)\bigr)\rightarrow\Bigl(\Box\bigl((\phi\wedge\eta)\rightarrow(\psi\wedge\eta)\bigr)\rightarrow\bigl(\Box\neg\eta\vee\Box(\phi\rightarrow\psi)\bigr)\Bigr)\;.
$$
Interpretation: unless it is empty, a proposition may be removed from a logical equation, when it appears in the both sides and is independent with the equation components.
\begin{description}
\item[Proof.]
It is easy to prove $(\phi\wedge\eta)\rightarrow(\psi\wedge\eta)\equiv\neg\eta\vee(\phi\rightarrow\psi)$\,.\\
Then  $\vdash\Box\bigl((\phi\wedge\eta)\rightarrow(\psi\wedge\eta)\bigr)\rightarrow\Box\bigl(\neg\eta\vee(\phi\rightarrow\psi)\bigr)$\,.\\
Now $\vdash\bigl((\phi\times\eta)\wedge(\psi\times\eta)\bigr)\rightarrow\bigl((\phi\rightarrow\psi)\times\neg\eta\bigr)$\,,
by independence invariance and b6.weak.A.\\
The proof is achieved by means of the preceding property, \emph{independence and proof}.
\item[$\Box\Box\Box$]\rien
\end{description}
\emph{Corollary.} $\vdash\phi\times\eta$\,, $\vdash\psi\times\eta$\,, $\vdash\Diamond\eta$ and $\phi\wedge\eta\equiv\psi\wedge\eta$ implies $\phi\equiv\psi$\,.
\\[5pt]
\emph{Corollary 2.} Being given $\psi$ and $\phi$ such that $\vdash\Diamond\phi$, the proposition $(\psi|\phi)$ is uniquely defined as the solution of equation $X\wedge\phi\equiv\psi\wedge\phi$ (with unknown $X$) which is independent of $\phi$.
\\[5pt]
This uniqueness is fake however, since the definition of $\times$ depends on $(|)$.
\begin{description}
\item[Proof.]
Let hypothesize $\vdash X\times\phi$ and $\vdash\Diamond\phi$.\\
Since $\vdash(\psi|\phi)\times\phi$ and $\psi\wedge\phi\equiv(\psi|\phi)\wedge\phi$, it is deduced from $X\wedge\phi\equiv\psi\wedge\phi$ that $X\equiv(\psi|\phi)$.
\item[$\Box\Box\Box$]\rien
\end{description}
\paragraph{Propositions equivalences [b6].}
$\vdash\Box(\psi\leftrightarrow\eta)\rightarrow\Box\bigl((\phi|\psi)\leftrightarrow(\phi|\eta)\bigr)$ (proved with b6 but without b6.weak.B).
\\[5pt]
Interpretation: equivalence is compliant with the conditioning.
\begin{description}
\item[Proof.]
First notice that all previous properties are obtained without b6.weak.B, and are thus obtained from b6.
\\[5pt]
From $(\phi|\psi)\wedge\psi\equiv\phi\wedge\psi$, $(\phi|\eta)\wedge\eta\equiv\phi\wedge\eta$ and $\vdash(\psi\leftrightarrow\eta)\rightarrow\bigl((\phi\wedge\psi)\leftrightarrow(\phi\wedge\eta)\bigr)$, it is deduced $\vdash(\psi\leftrightarrow\eta)\rightarrow\Bigl(\bigl((\phi|\psi)\wedge\psi\bigr)\leftrightarrow\bigl((\phi|\eta)\wedge\eta\bigr)\Bigr)\,.$\\
Then $\vdash(\psi\leftrightarrow\eta)\rightarrow\Bigl(\bigl((\phi|\psi)\wedge\psi\bigr)\leftrightarrow\bigl((\phi|\eta)\wedge\psi\bigr)\Bigr)$ and finally:\vspace{-5pt}
$$
\vdash\Box(\psi\leftrightarrow\eta)\rightarrow\Box\Bigl(\bigl((\phi|\psi)\wedge\psi\bigr)\leftrightarrow\bigl((\phi|\eta)\wedge\psi\bigr)\Bigr)\,.
$$
Now $\vdash(\phi|\psi)\times\psi$ and $\vdash(\phi|\eta)\times\eta$\,;
since $\vdash\Box(\psi\leftrightarrow\eta)\rightarrow\Bigr(\bigl(\eta\times(\phi|\eta)\bigr)\leftrightarrow\bigl(\psi\times(\phi|\eta)\bigr)\Bigl)$, it comes $\vdash\Box(\psi\leftrightarrow\eta)\rightarrow\bigl((\phi|\eta)\times\psi\bigr)$ by b6.
\\
Finally $\vdash\Box(\psi\leftrightarrow\eta)\rightarrow\biggl(\bigl((\phi|\psi)\times\psi\bigr)\wedge\bigl((\phi|\eta)\times\psi\bigr)\wedge\Box\Bigl(\bigl((\phi|\psi)\wedge\psi\bigr)\leftrightarrow\bigl((\phi|\eta)\wedge\psi\bigr)\Bigr)\biggr)$\,.\\
Applying the regularity, it comes $\vdash\Box(\psi\leftrightarrow\eta)\rightarrow\Bigl(\Box\neg\psi\vee\Box\bigl((\phi|\psi)\leftrightarrow(\phi|\eta)\bigr)\Bigr)$\,.\\[5pt]
Now $\vdash\Box(\psi\leftrightarrow\eta)\rightarrow (\Box\neg\psi\leftrightarrow\Box\neg\eta)$
and $\vdash\Box\neg\Xi\rightarrow\Box\bigl((\phi|\Xi)\leftrightarrow\phi\bigr)$ for $\Xi=\psi$ or $\eta$.\\
It is deduced $\vdash\bigl(\Box\neg\psi\wedge\Box(\psi\leftrightarrow\eta)\bigr)\rightarrow\Box\bigl((\phi|\psi)\leftrightarrow(\phi|\eta)\bigr)$\,, thus completing the proof.
\item[$\Box\Box\Box$]\rien
\end{description}
\emph{Corollary.} Axiom b6 implies b6.weak.B.
In particular, DmBL$_\ast$ is weaker than DmBL.
\\[5pt] 
\emph{Corollary of b6 or b6.weak.B.} $\psi\equiv\eta$ implies $(\phi|\psi)\equiv(\phi|\eta)$.
\\[5pt]
Together with the properties of the system $T$ and \emph{left equivalences}, this last result implies that the equivalence relation $\equiv$ is compliant with the logical operators of DmBL/DmBL$_\ast$.
In particular, replacing a sub-proposition with an equivalent sub-proposition within a theorem still makes a theorem.
\paragraph{Reduction rule [b6].}
The axiom b6 implies $\bigl(\phi\big|(\psi|\phi)\bigr)\equiv\phi$\,.
\begin{description}
\item[Proof.]
Since $\vdash(\psi|\phi)\times\phi$, b6 implies $\vdash\phi\times(\psi|\phi)$ and $\vdash \Box\Bigl(\bigl(\phi\big|(\psi|\phi)\bigr)\leftrightarrow\phi\Bigr)$\,.
\item[$\Box\Box\Box$]\rien
\end{description}
\paragraph{Markov Property [b6].}
$$
\vdash\left(\left(
\bigwedge_{\tau=1}^{t-2}\bigl((\phi_t|\phi_{t-1})\times\phi_{\tau}\bigr)
\right)\wedge\Diamond\left(\bigwedge_{\tau=1}^{t-1}\phi_{\tau}\right)\right)\longrightarrow\Box\left(
(\phi_t|\phi_{t-1})\leftrightarrow\left(\phi_t\left|\bigwedge_{\tau=1}^{t-1}\phi_{\tau}\right.\right)
\right)\;.
$$
Interpretation: the Markov property holds, when the conditioning is independent of the past and the past is possible.
\begin{description}
\item[Proof.]
Since $\vdash(\phi_t|\phi_{t-1})\times\phi_{t-1}$, it comes $\vdash\left(
\bigwedge_{\tau=1}^{t-2}\bigl((\phi_t|\phi_{t-1})\times\phi_{\tau}\bigr)
\right)\rightarrow\left(
\bigwedge_{\tau=1}^{t-1}\bigl((\phi_t|\phi_{t-1})\times\phi_{\tau}\bigr)
\right)$\,.\\
Then $\vdash\left(
\bigwedge_{\tau=1}^{t-2}\bigl((\phi_t|\phi_{t-1})\times\phi_{\tau}\bigr)
\right)\rightarrow\left((\phi_t|\phi_{t-1})\times \left(\bigwedge_{\tau=1}^{t-1}\phi_\tau\right)\right)$\,.\\
Now, $(\phi_t|\phi_{t-1})\wedge \left(\bigwedge_{\tau=1}^{t-1}\phi_\tau\right)\equiv\bigwedge_{\tau=1}^{t}\phi_\tau
\equiv\left(\phi_t\left|\bigwedge_{\tau=1}^{t-1}\phi_\tau\right.\right)\wedge \left(\bigwedge_{\tau=1}^{t-1}\phi_\tau\right)$\,.\\
Since $\vdash\left(\phi_t\left|\bigwedge_{\tau=1}^{t-1}\phi_\tau\right.\right)\times \left(\bigwedge_{\tau=1}^{t-1}\phi_\tau\right)$\,, the proof is achieved by applying the regularity.
\item[$\Box\Box\Box$]\rien
\end{description}
%
%
\paragraph{Link between $\bigl((\eta|\psi)\big|\phi\bigr)$ and $(\eta|\phi\wedge\psi)$ [b6]}.\rien\\
\emph{Is proved} $\bigl((\eta|\psi)\big|\phi\bigr)\wedge(\phi\wedge\psi)\equiv(\eta|\phi\wedge\psi)\wedge(\phi\wedge\psi)\,.$
\begin{description}
\item[Proof.]
$$
\bigl((\eta|\psi)\big|\phi\bigr)\wedge\phi\wedge\psi\equiv(\eta|\psi)\wedge\phi\wedge\psi\equiv\phi\wedge\psi\wedge\eta\equiv(\eta|\phi\wedge\psi)\wedge(\phi\wedge\psi)\;.
$$
\item[$\Box\Box\Box$]\rien
\end{description}
This is a quite limited result, and it is \emph{tempting} to assume the additional but \emph{controversial} axiom ``$\bigl((\eta|\psi)\big|\phi\bigr)\equiv(\eta|\phi\wedge\psi)\quad\mbox{\small$(\ast)$}$''\,.
There is a really critical point here.
A main argument against \emph{Bayesian propositions} relies on \emph{Lewis' negative result}\,.
This result will be explained later, when speaking about probabilities.
It will be shown that \emph{Lewis' negative result} does not hold with DmBL.
However this becomes false with~$(\ast)$, and actually axiom~$(\ast)$ implies a strictly logical counterpart to \emph{Lewis}\,:
\begin{quote}
\emph{Let $\bigl((\eta|\psi)\big|\phi\bigr)\equiv(\eta|\phi\wedge\psi)\quad\mbox{\small$(\ast)$}$ be assumed as an axiom.\\
Then $\vdash\Diamond(\phi\wedge\psi)\rightarrow
\bigl(\Box(\phi\leftrightarrow\psi)\vee(\phi\times\psi)\bigr)$\,.}
\end{quote}
Interpretation: if $\phi$ and $\psi$ are not exclusive and not equivalent, then they are independent.
This is irrelevant and forbids the use of axiom~$(\ast)$.
\begin{description}
\item[Proof.]
Since $\neg(\phi\rightarrow\psi)\equiv\phi\wedge\neg\psi$\,, it is equivalent to prove:
$$
\vdash\Bigl(\Diamond(\phi\wedge\psi)\wedge\bigl(\Diamond(\neg\psi\wedge\phi)\vee\Diamond(\neg\phi\wedge\psi)\bigr)\Bigr)\rightarrow(\phi\times\psi)
\;.
$$
Since $\times$ is symmetric, it is sufficient to prove $\vdash\bigl(\Diamond(\phi\wedge\psi)\wedge\Diamond(\neg\psi\wedge\phi)\bigr)\rightarrow(\phi\times\psi)$\,.\\[5pt]
The \emph{introspection} property implies $\vdash\Diamond(\phi\wedge\psi)\rightarrow\Box(\phi\wedge\psi|\phi\wedge\psi)$, denoted $(a)$, and $\vdash\Diamond(\neg\psi\wedge\phi)\rightarrow\Box(\neg\psi\wedge\phi|\neg\psi\wedge\phi)$\,. \\
It is thus deduced $\vdash\Diamond(\phi\wedge\psi)\rightarrow\Box\Bigr((\psi|\psi\wedge\phi)\leftrightarrow\bigr((\psi|\psi\wedge\phi)\wedge(\psi\wedge\phi|\psi\wedge\phi)\bigl)\Bigl)$, denoted $(b)$, and $\vdash\Diamond(\neg\psi\wedge\phi)\rightarrow\Box\Bigr((\psi|\neg\psi\wedge\phi)\leftrightarrow\bigr((\psi|\neg\psi\wedge\phi)\wedge(\neg\psi\wedge\phi|\neg\psi\wedge\phi)\bigl)\Bigl)$. \\
From the deduction $(\psi|\neg\psi\wedge\phi)\wedge(\neg\psi\wedge\phi|\neg\psi\wedge\phi)\equiv(\bot|\neg\psi\wedge\phi)\equiv\bot$, it is derived $\vdash\Diamond(\neg\psi\wedge\phi)\rightarrow\Box\bigl((\psi|\neg\psi\wedge\phi)\leftrightarrow\bot\bigr)$, denoted $(c)$. \\
From the deduction $(\psi|\psi\wedge\phi)\wedge(\psi\wedge\phi|\psi\wedge\phi)\equiv(\psi\wedge\phi|\psi\wedge\phi)$, $(a)$ and $(b)$, it comes $\vdash\Diamond(\phi\wedge\psi)\rightarrow\Box\bigl((\psi|\psi\wedge\phi))\leftrightarrow\top\bigr)$, denoted $(d)$.\\
Now $(\psi|\phi)\equiv\bigl((\psi|\phi)\wedge\psi\bigr)\vee\bigl((\psi|\phi)\wedge\neg\psi\bigr)\equiv
\Bigl(\bigl((\psi|\phi)\big|\psi\bigr)\wedge\psi\Bigr)\vee\Bigl(\bigl((\psi|\phi)\big|\neg\psi\bigr)\wedge\neg\psi\Bigr)$\,, and by applying axiom $(\ast)$, $(\psi|\phi)\equiv\bigl((\psi|\phi\wedge\psi)\wedge\psi\bigr)\vee\bigl((\psi|\phi\wedge\neg\psi)\wedge\neg\psi\bigr)$.\\
Then $\vdash\bigl(\Diamond(\neg\psi\wedge\phi)\wedge\Diamond(\phi\wedge\psi)\bigr)\rightarrow\Box\Bigl(
(\psi|\phi)\leftrightarrow\bigl((\top\wedge\psi)\vee(\bot\wedge\neg\psi)\bigr)
\Bigr)$ by $(c)$ and $(d)$. \\
At last $\vdash\bigl(\Diamond(\neg\psi\wedge\phi)\wedge\Diamond(\phi\wedge\psi)\bigr)\rightarrow\Box\bigl(
(\psi|\phi)\leftrightarrow\psi\bigr)$\,.
\item[$\Box\Box\Box$]\rien
\end{description}
\paragraph{Some comparisons with CEA.}
It is interesting to compare our logical results with what is obtained in Conditional Event Algebra.
There are different possible CEA, but we will restrict our discussion to the DGNW and later to the Product Space CEA \cite{goodman}.
In appendix~\ref{BayesMod:Comp&Ex} it will also be discussed about the Product Space CEA \cite{goodman,goodman2}, by a comparison with the conditional model of DmBL$_\ast$.
\\[3pt]
The classical nature of the sub-universes are a common value for the CEAs and DmBL.
Now, DGNW introduces some new equivalences, which are not in DmBL: 
\begin{itemize}
\item $(a|b)\wedge(c|d)=\bigl(a\wedge b\wedge c\wedge d\big|(\neg a\wedge b)\vee(\neg c\wedge d)\vee (b\wedge d)\bigl)$\,,
\item $(a|b)\vee(c|d)=\bigl((a\wedge b)\vee(c\wedge d)\big|(a\wedge b)\vee(c\wedge d)\vee (b\wedge d)\bigr)$\,,
\end{itemize}
Owing to the free model constructed in next section, it is likely that these equivalences will need some additional axioms in DmBL.
\\[3pt]
Another interesting difference is noticed for the conditioning.
While DGNW provides:
$$
(a|b\wedge c)\wedge (b|c)=(a\wedge b|c)\quad\mbox{and}\quad(a|b)\wedge (b|\Omega)=(a\wedge b|\Omega)\,,
$$
DmBL implies $(\psi|\phi)\wedge \phi\equiv\phi\wedge \psi$\,, \emph{but not} $(\eta|\phi\wedge\psi)\wedge (\psi|\phi)\equiv(\psi\wedge \eta|\phi)$\,.
\\[3pt]
Besides, DmBL allows the trivial conditioning to be removed, \emph{i.e.} $(\phi|\top)\equiv\phi$, while $(a|\Omega)=a$ is not allowed in DGNW.
About the Product Space CEA, the answer is less clear since \cite{goodman2} removes it, while \cite{goodman} does not.
In will be discussed more precisely about this question in section~\ref{Section:Proba:DmBL}, when considering Lewis's triviality.
\\[3pt]
A last point is about the closure of DGNW.
It is proposed in \cite{goodman}: 
$$
\bigl((a|b)\big|(c|d)\bigr)=\big(a\big|(b\wedge\neg a\wedge\neg d)\vee(b\wedge c\wedge d)\bigr)\;.
$$
But this definition fails to satisfy the intuitive relation:
$$
P\bigl((a|b)\big|(c|d)\bigr)=P\bigl((a|b)\wedge(c|d)\bigr)/P(c|d)\;.
$$
It is proved in section~\ref{Section:Proba:DmBL} that the probabilistic extension of DmBL verifies this property.
\section{Models}
\label{Section:Model:DmBL}
\subsection{Toward a Model}
In this paragraph, it is discussed about the link between Kripke models for DmBL/DmBL$_\ast$ and a more basic structure called \emph{conditional models}.
\paragraph{Definition.}
A Kripke model for DmBL (respectively DmBL$_\ast$) is a quadruplet $(W,R,H,f)$ such that $W$ is a set of worlds, $R\subset W\times W$ is the accessibility relation, $H:\mathcal{L}\longrightarrow\mathcal{P}(W)$ is an assignment function, $f:H(\mathcal{L})\times H(\mathcal{L})\longrightarrow\mathcal{P}(W)$, and verifying:
\begin{itemize}
\item $H(\phi\wedge\psi)=H(\phi)\cap H(\psi)$\,,
\item $H(\neg\phi)=W\setminus H(\phi)$\,,
\item $H(\top)=W$\,,
\item $H(\phi\vee\psi)=H(\phi)\cup H(\psi)$ and $H(\phi\rightarrow\psi)=\bigl(W\setminus H(\phi)\bigr)\cup H(\psi)$\,,
\item $H(\Box\phi)=\bigl\{t\in W\big/ \forall u\in W,\,(t,u)\in R\Rightarrow u\in H(\phi)\bigr\}$\,,
\item $H\bigl((\psi|\phi)\bigr)=f\bigl(H(\psi),H(\phi)\bigr)$\,,
\item $H(\phi)=W$ for any $\phi$ such that $\vdash\phi$ is an axiom of the form m3, b1, b2, b3, b4,\\
or b6 (respectively b6.weak.A and b6.weak.B).
\end{itemize}
Notice that the rules and axioms c$\ast$, modus ponens, m1 and m2 and are compliant with the model by construction.
m4 and b5 are just definitions.
\paragraph{Definition 2.} A conditional model for DmBL (respectively DmBL$_\ast$) is a quadruplet $(W,M,h,f)$ such that $W$ is a set of worlds, $M\subset\mathcal{P}(W)$ is a set of admissible propositions, $h:\Theta\longrightarrow M$ is an assignment function, $f:M\times M\longrightarrow M$, and verifying:
\begin{description}
\item[$\rien\quad\bullet$] $M$ is a Boolean sub-algebra of $\mathcal{P}(W)$, \emph{i.e.} $A\cap B\in M$ and $W\setminus A\in M$ for any $A,B\in M$\,,
\item[$\rien\quad\beta1$.] $A\subset B$ and $A\ne\emptyset$ imply $f(B,A)=W$\,, for any $A,B\in M$\,,
\item[$\rien\quad\beta2$.] $f(B\cup C,A)\subset f(B,A)\cup f(C,A)$\,, for any $A,B,C\in M$\,,
\item[$\rien\quad\beta3$.] $A\cap f(B,A) \subset B$\,, for any $A,B\in M$\,,
\item[$\rien\quad\beta4$.] $f(W\setminus B,A)=W\setminus f(B,A)$\,,
for any $A,B\in M$\,,
\item[$\rien\quad\beta6$.] $f(B,A)=B$ implies $f(A,B)=A$\\
(respectively $\beta6w.$ $f(B,A)=B$ implies  $f(B,W\setminus A)=B$)\,,
for any $A,B\in M$\,.
\end{description}
Such model does not implement the modalities.
\\[5pt]
Remark: in both models, the function $f$ is the representation of the conditional $(|)$\,.
\paragraph{Model transfer.}
Let $(W,M,h,f)$ be a conditional model for DmBL (respectively DmBL$_\ast$).
Let $R=W\times W$ and define $H$ by:
\begin{itemize}
\item $H(\theta)=h(\theta)$ for any $\theta\in\Theta$\,,
\item $H(\phi\wedge\psi)=H(\phi)\cap H(\psi)$\,,
$H(\neg\phi)=W\setminus H(\phi)$\,,
$H(\top)=W$\,,
\item $H(\phi\vee\psi)=H(\phi)\cup H(\psi)$ and $H(\phi\rightarrow\psi)=\bigl(W\setminus H(\phi)\bigr)\cup H(\psi)$\,,
\item $H(\phi)=W\Rightarrow H(\Box\phi)=W$ and $H(\phi)\ne W\Rightarrow H(\Box\phi)=\emptyset$\,,
\item $H\bigl((\psi|\phi)\bigr)=f\bigl(H(\psi),H(\phi)\bigr)$\,.
\end{itemize}
Then $(W,R,H,f)$ is a Kripke model for DmBL (respectively DmBL$_\ast$).
\begin{description}
\item[Proof.] First notice that the above construction of $H$ is possible for any proposition $\phi\in\mathcal{L}$, since it is always obtained $H(\phi)\in M$\,.\\
Now, $R=W\times W$ implies $H(\Box\phi)=\bigl\{t\in W\big/ \forall u\in W,\,(t,u)\in R\Rightarrow u\in H(\phi)\bigr\}$\,, so that $(W,R,H,f)$ is actually a Kripke model.\\
Let verify the compliance with m3, b1, b2, b3, b4 and b6 (resp. b6.weak.$\ast$).\\[10pt]
\underline{Compliance with m3} is obtained from the fact that $R$ is reflexive.
\\[5pt]
Proof of $H\Bigl(\Box(\phi\rightarrow\psi)\rightarrow\bigl(\Box\neg\phi\vee\Box(\psi|\phi)\bigr)\Bigr)=W$\,, \emph{i.e.} \underline{compliance with b1}.\\
By definition, $H\Bigl(\Box(\phi\rightarrow\psi)\rightarrow\bigl(\Box\neg\phi\vee\Box(\psi|\phi)\bigr)\Bigr)=$\\
\rien\hspace{150pt}$H\bigl(\Box(\psi|\phi)\bigr)\cup H(\Box\neg\phi) \cup \Bigl(W\setminus H\bigl(\Box(\phi\rightarrow\psi)\bigr)\Bigl)$\,.\\
Then $H\bigl(\Box(\phi\rightarrow\psi)\bigr)=\emptyset$ or $H(\Box\neg\phi)=W$ imply \mbox{$H\Bigl(\Box(\phi\rightarrow\psi)\rightarrow\bigl(\Box\neg\phi\vee\Box(\psi|\phi)\bigr)\Bigr)=W$\,.\hspace{-20pt}\rien}
Otherwise $H\bigl(\Box(\phi\rightarrow\psi)\bigr)\ne\emptyset$ and $H(\Box\neg\phi)\ne W$\,, thus implying $H(\phi\rightarrow\psi)=W$ and $H(\neg\phi)\ne W$\,.\\
It is deduced $H(\phi)\cap \bigl(W\setminus H(\psi)\bigr)=\emptyset$ and $H(\phi)\ne\emptyset$\,.\\
Then $H(\phi)\subset H(\psi)$ and $H(\phi)\ne\emptyset$\,, and $f\bigl(H(\psi), H(\phi)\bigr)=W$ by using $\beta1$.\\
Finally $H\bigl(\Box(\psi|\phi)\bigr)=W$ and $H\Bigl(\Box(\phi\rightarrow\psi)\rightarrow\bigl(\Box\neg\phi\vee\Box(\psi|\phi)\bigr)\Bigr)=W$\,.
\\[5pt]
Proof of $H\bigl(\neg(\neg\psi|\phi)\leftrightarrow(\psi|\phi)\bigr)=W$\,, \emph{i.e.} \underline{compliance with b4}.\\
It is deduced $
H\bigl((\neg\psi|\phi)\bigr)=f\bigl(H(\neg\psi),H(\phi)\bigr)=f\bigl(W\setminus H(\psi),H(\phi)\bigr)=W\setminus f\bigl(H(\psi),H(\phi)\bigr)=W\setminus H\bigl((\psi|\phi)\bigr)
$\,, by using $\beta4$.\\
Then $\Bigl(H\bigl((\neg\psi|\phi)\bigr)\cap \bigl(W\setminus H\bigl((\psi|\phi)\bigr)\bigr)\Bigr)\cup\Bigl(\bigl(W\setminus H\bigl((\neg\psi|\phi)\bigr)\bigr)\cap H\bigl((\psi|\phi)\bigr)\Bigr)=W$\,.\\
Then $H\bigl(\neg(\neg\psi|\phi)\leftrightarrow(\psi|\phi)\bigr)=W$\,.
\\[5pt]
Proof of $H\Bigl((\psi\rightarrow\eta|\phi)\rightarrow\bigl((\psi|\phi)\rightarrow(\eta|\phi)\bigr)\Bigr)=W$\,, \emph{i.e.} \underline{compliance with b2}.\\
It is deduced $H\Bigl((\psi\rightarrow\eta|\phi)\rightarrow\bigl((\psi|\phi)\rightarrow(\eta|\phi)\bigr)\Bigr)=
H\Bigl(\neg(\neg\psi\vee\eta|\phi)\vee\bigl(\neg(\psi|\phi)\vee(\eta|\phi)\bigr)\Bigr)$\\
\rien\hspace{50pt}$=\Bigl(W\setminus\Bigl(H\bigl((\neg\psi\vee\eta|\phi)\bigr)\cap H\bigl((\psi|\phi)\bigr)\Bigr)\cup H\bigl((\eta|\phi)\bigr)$\,.\\
Now $H\bigl((\neg\psi\vee\eta|\phi)\bigr)=f\bigl(H(\neg\psi\vee\eta),H(\phi)\bigr)=f\bigl(H(\neg\psi)\cup H(\eta),H(\phi)\bigr)=f\bigl(H(\neg\psi),H(\phi)\bigr)$\\
\rien\hspace{25pt}$\cup f\bigl(H(\eta),H(\phi)\bigr)=H\bigl((\neg\psi|\phi)\bigr)\cup H\bigl((\eta|\phi)\bigr) =\Bigl(W\setminus H\bigl((\psi|\phi)\bigr)\Bigr)\cup H\bigl((\eta|\phi)\bigr)$\,, by $\beta2$.\\
Then \mbox{$H\Bigl((\psi\rightarrow\eta|\phi)\rightarrow\bigl((\psi|\phi)\rightarrow(\eta|\phi)\bigr)\Bigr)=\Bigl(W\setminus\Bigl(H\bigl((\psi|\phi)\bigr)\cap H\bigl((\eta|\phi)\bigr)\Bigr)\Bigr)\cup H\bigl((\eta|\phi)\bigr)=W\,.$\hspace{-20pt}\rien}
\\[5pt]
Proof of $H\bigl((\psi|\phi)\rightarrow(\phi\rightarrow\psi)\bigr)=W$\,, \emph{i.e.} \underline{compliance with b3}.\\
Immediate from $\beta3$, \emph{i.e.} $H(\phi)\cap f\bigl(H(\psi),H(\phi)\bigr)\subset H(\psi)$\,.
\\[5pt]
{\bf Case DmBL.}
Proof of $H\Bigl(\Box\bigl((\psi|\phi)\leftrightarrow\psi\bigr)\leftrightarrow\Box\bigl((\phi|\psi)\leftrightarrow\phi\bigr)\Bigr)=W$\,, \emph{i.e.} \underline{compliance with b6}.\\
By $\beta6$, $f\bigl(H(\psi),H(\phi)\bigr)=H(\psi)$ if and only if $f\bigl(H(\phi),H(\psi)\bigr)=H(\phi)$\,.\\
Then $H\bigl((\psi|\phi)\leftrightarrow\psi\bigr)=W$ if and only if $H\bigl((\phi|\psi)\leftrightarrow\phi\bigr)=W$\,.\\
As a consequence, $H\Bigl(\Box\bigl((\psi|\phi)\leftrightarrow\psi\bigr)\Bigr)=H\Bigl(\Box\bigl((\phi|\psi)\leftrightarrow\phi\bigr)\Bigr)$ and the result.
\\[5pt]
{\bf Case DmBL$_\ast$.}\\
Proof of $H\Bigl(\Box\bigl((\psi|\phi)\leftrightarrow\psi\bigr)\leftrightarrow\Box\bigl((\psi|\neg\phi)\leftrightarrow\psi\bigr)\Bigr)=W$\,, \emph{i.e.} \underline{compliance with b6.weak.A}.\\
By $\beta6w$, $f\bigl(H(\psi),H(\phi)\bigr)=H(\psi)$ if and only if $f\bigl(H(\psi),W\setminus H(\phi)\bigr)=H(\psi)$\,.\\
Then $H\bigl((\psi|\phi)\leftrightarrow\psi\bigr)=W$ if and only if $H\bigl((\psi|\neg\phi)\leftrightarrow\psi\bigr)=W$\,.\\
As a consequence, $H\Bigl(\Box\bigl((\psi|\phi)\leftrightarrow\psi\bigr)\Bigr)=H\Bigl(\Box\bigl((\psi|\neg\phi)\leftrightarrow\psi\bigr)\Bigr)$ and the result.
\\[5pt]
Proof of $H\Bigl(\Box(\psi\leftrightarrow\eta)\rightarrow\Box\bigl((\phi|\psi)\leftrightarrow(\phi|\eta)\bigr)\Bigr)=W$\,, \emph{i.e.} \underline{compliance with b6.weak.B}.\\
Assume first $H(\psi)\ne H(\eta)$\,.\\
It comes $H(\psi\leftrightarrow\eta)\ne W$ and $H\bigl(\Box(\psi\leftrightarrow\eta)\bigr)=\emptyset$\,.\\
Then $H\Bigl(\Box(\psi\leftrightarrow\eta)\rightarrow\Box\bigl((\phi|\psi)\leftrightarrow(\phi|\eta)\bigr)\Bigr)=W$\,.\\[3pt]
Assume now $H(\psi)= H(\eta)$\,.\\
Then $H\bigl((\phi|\psi)\bigr)=H\bigl((\phi|\eta)\bigr)$ and $H\Bigl(\Box\bigl((\phi|\psi)\leftrightarrow(\phi|\eta)\bigr)\Bigr)=H\bigl((\phi|\psi)\leftrightarrow(\phi|\eta)\bigr)=W$\,.\\
At last $H\Bigl(\Box(\psi\leftrightarrow\eta)\rightarrow\Box\bigl((\phi|\psi)\leftrightarrow(\phi|\eta)\bigr)\Bigr)=W$\,.
\item[$\Box\Box\Box$]\rien
\end{description}
Conditional models are defined from the classical and conditional operators only.
In fact, they have been set first for a non modal construction of the Bayesian logic\cite{Dambreville:DBL}, and were complete in this case.
It appears that conditional models are translated into DmBL/DmBL$_\ast$ Kripke models.
The derived models are of course not complete for DmBL/DmBL$_\ast$, but are sufficient to prove that the conditional operator cannot be reduced to trivialities.
\\[5pt]
Moreover, the \emph{model transfer property} also suggests that the conditional operator is \emph{not} constructed from the modal operator: it is even possible to construct $(|)$ when $\Box$ is trivial in the model ($R=W\times W$ implies $H(\phi)=W$ or $H(\Box\phi)=\emptyset$)\,.
\\[5pt]
In this paper a free conditional model is constructed for DmBL$_\ast$\,, with some completeness results.
In section~\ref{Section:Proba:DmBL} a model for DmBL will be also derived but not constructed.
This model of DmBL will be non trivial, but no completeness result will be derived.
\subsection{Construction of a free conditional model for DmBL$_\ast$}
In the sequel, $\Theta$ is assumed to be finite.
A free conditional model for DmBL$_\ast$ will be constructed as a limit of partial models.
These models are constructed recursively, based on the iteration of $(|)$ on any propositions.
\subsubsection{Definition of partial models}
In this section are constructed a sequence $(\Omega_n, M_n, h_n, f_n,\Lambda_n)_{n\in\Nset}$ and a sequence of one-to-one morphisms $(\mu_n)_{n\in\Nset}$ such that:
\begin{itemize}
\item $M_n$ is a Boolean sub-algebra of $\mathcal{P}(\Omega_n)$, $h_n:\Theta\rightarrow M_n$ and $f_n:M_n\times M_n\rightarrow M_n$\,,\\
($f_n$ will be partially defined)
\item $\mu_n:M_n\rightarrow M_{n+1}$ is such that $\mu_n(A\cap B)=\mu_n(A)\cap\mu_n(B)$, $\mu_n(\Omega_n\setminus A)=\Omega_{n+1}\setminus\mu_n(A)$ (\emph{i.e.} $\mu_n$ is a Boolean morphism) and $\forall\theta\in\Theta,\,h_{n+1}(\theta)=\mu_n\bigl(h_n(\theta)\bigr)$\,,
\item For any $A,B\in M_n$ such that $f_n(B,A)$ is defined, then $f_{n+1}\bigl(\mu_n(B),\mu_n(A)\bigr)$ is defined and $f_{n+1}\bigl(\mu_n(B),\mu_n(A)\bigr)=\mu_n\bigl(f_n(B,A)\bigr)$\,,
\item $\Lambda_n$ is a list of elements of $M_n$, which is used as a task list of the construction (refer to the subsequent paragraphs.).
\end{itemize}
\emph{Remark.}
The functions $f_n$ represent the partial construction of $(|)$\,.
The morphisms $\mu_n$ characterize the ``inclusion'' of the partial models.
\\[5pt]
In a subsequent section, a conditional model $\bigl(\Omega_\infty,M_\infty,h_\infty,f_\infty\bigr)$ will be defined as the limit of $(\Omega_n, M_n, h_n, f_n)_{n\in\Nset}$ associated to $(\mu_n)_{n\in\Nset}$\,.
\paragraph{Notations and definitions.}
For any $A\in M_n$, it is defined $\sim A=\Omega_n\setminus A$\,.\\[5pt]
Any singleton $\{\omega\}$ may be denoted $\omega$ if the context is not ambiguous.\\[5pt]
For any $m>n$ and $A\in M_n$\,, it is defined $A_{[m]}=\mu_{m-1}\circ\dots\circ\mu_{n}(A)$\,.\\[5pt]
The Cartesian product of sets $A$ and $B$ is denoted $A\times B$\,;
the functions $\mathrm{id}$ and $T$ are defined over pairs by $\mathrm{id}(x,y)=(x,y)$ and $T(x,y)=(y,x)$\,.
\paragraph{Initialization.}
Define $(\Omega_0, M_0, h_0, f_0,\Lambda_0)$ by:
\begin{itemize}
\item $\Omega_0=\{0,1\}^{\Theta}$,
\item $M_0=\mathcal{P}(\Omega_0)$,
\item $h_0(\theta)=\bigl\{(\delta_\tau)_{\tau\in\Theta}\in\Omega_0\,\big/\,\delta_\theta=1\bigr\}$\,,
\item $f_0(A,\emptyset)=f_0(A,\Omega_0)=A$ for any $A\in M_0$\,,
\item $\Lambda_0=(s_0,f_0,\lambda_0)$ is a list defined by $s_0=0$, $f_0=\mathrm{card}(M_0)-2$ and $\lambda_0$ is a one-to-one mapping from $[\![s_0,f_0-1]\!]$ to $M_0\setminus\{\emptyset,\Omega_0\}$\,,
such that $\lambda_0(2t)=\sim\lambda_0(2t+1)\,,\;\forall t$\,.
\end{itemize}
\paragraph{Step $n$ to step $n+1$.}\rien\\
Let $(\Omega_k, M_k, h_k, f_k,\Lambda_k)_{0\le k\le n}$ and the one-to-one morphisms $(\mu_k)_{0\le k\le n-1}$ be constructed.\\[5pt]
Notice that $\lambda_n(s_n)=\sim\lambda_n(s_n+1)$ by construction of $\Lambda_n$\,.\\
\underline{Define $b_{n}=\lambda_n(s_n)$}\,.\\
Then, construct the set $I_n$ and the sequences $\Gamma_n(i),\Pi_n(i)|_{i\in I_n}$ according to the cases:
\subparagraph{Case 0.} There is $m<n$ such that $\{b_{m[n]},\sim b_{m[n]}\}=\{b_n,\sim b_n\}$\,.\\
Then, notice that $b_n=b_{m[n]}$ by construction of $\Lambda$\,.\\
Let $\nu$ be the greatest of such $m$;
then define $I_n=\mu_\nu(b_{\nu})\times\sim \mu_\nu(b_{\nu})$\,,\\
$\Pi_n(\omega,\omega')=f_n(\omega'_{[n]},\sim b_n)\cap\omega_{[n]}$ and $\Gamma_n(\omega,\omega')=f_n(\omega_{[n]},b_n)\cap\omega'_{[n]}$ for any $(\omega,\omega')\in I_n$\,.$^\dagger$
\\[5pt]
$^\dagger$Remark: case 0 means that the construction of $f(\cdot,b_n)$ and of $f(\cdot,\sim b_n)$ has already begun over the propositions of $M_{\nu+1}$.
\subparagraph{Case 1.} Case 0 does not hold;\\
Define $I_n=\{b_n\}$\,, $\Pi_n(i)=i$ 
and $\Gamma_n(i)=\sim i$ for any $i\in I_n$\,.
\\[5pt]
Remark: case 1 means that $f(\cdot,b_n)$ and $f(\cdot,\sim b_n)$ are constructed for the first time.
\subparagraph{Setting.}
$(\Omega_{n+1}, M_{n+1}, h_{n+1}, f_{n+1},\Lambda_{n+1})$ and $\mu_n$ are defined by:
\begin{itemize}
\item $\mu_n(A)=\bigcup_{i\in I_n} \biggl(\Bigl(\bigl(A\cap\Pi_n(i)\bigr)\times\Gamma_n(i)\Bigr)
\cup \Bigl(\bigl(A\cap\Gamma_n(i)\bigr)\times\Pi_n(i)\Bigr)\biggr)$ for any $A\in M_n$\,,
\item $\Omega_{n+1}=\mu_n(\Omega_n)$\,,
\item $\forall\theta\in\Theta,\,h_{n+1}(\theta)=\mu_n\bigl(h_n(\theta)\bigr)$\,,
\item $M_{n+1}=\mathcal{P}(\Omega_{n+1})$\,,
\item $f_{n+1}(A,\emptyset)=f_{n+1}(A,\Omega_{n+1})=A$ for any $A\in M_{n+1}$\,,
\item For any $A\in M_n\setminus\{b_n,\sim b_n,\emptyset,\Omega_n\}$ and any $B\in M_n$ such that $f_n(B,A)$ is defined, then $f_{n+1}\bigl(\mu_n(B),\mu_n(A)\bigr)$ is defined and $f_{n+1}\bigl(\mu_n(B),\mu_n(A)\bigr)=\mu_n\bigl(f_n(B,A)\bigr)$\,,
\item For any $A\in M_{n+1}$\,, set
$
f_{n+1}\bigl(A,\mu_n(b_n)\bigr)=(\mathrm{id}\cup T)
\biggl(A\cap\Bigl(\bigcup_{i\in I_n}\bigl(\Pi_n(i)\times\Gamma_n(i)\bigr)\Bigr)\biggr)
$\\[5pt]
and
$
f_{n+1}\bigl(A,\sim\mu_n(b_n)\bigr)=(\mathrm{id}\cup T)
\biggl(A\cap\Bigl(\bigcup_{i\in I_n}\bigl(\Gamma_n(i)\times\Pi_n(i)\bigr)\Bigr)\biggr)
\;,
$
\item $\Lambda_{n+1}=(s_{n+1},f_{n+1},\lambda_{n+1})$ is such that:
\begin{itemize}
\item $s_{n+1}=s_n+2=2n+2$, $f_{n+1}=s_{n+1}+\mathrm{card}(M_{n+1})-2$\,,
\item $\lambda_{n+1}$ is a one-to-one mapping from $[\![s_{n+1},f_{n+1}-1]\!]$ to $M_{n+1}\setminus\{\emptyset,\Omega_{n+1}\}$\,,
\item $\lambda_{n+1}(t)=\mu_n\bigl(\lambda_n(t)\bigr)$ for any $t\in[\![s_{n+1},f_{n}-1]\!]$\,,
\item $\lambda_{n+1}(2t)=\sim\lambda_{n+1}(2t+1)\,,\;\forall t$
and $\lambda_{n+1}(f_{n+1}-2)=\mu_n(b_n)$\,.
\end{itemize}
\emph{This definition ensures a cyclic and full construction of $f(\cdot,b_n)$ and $f(\cdot,\sim b_n)$}.\footnote{In regards to the mapping $\mu$, the list $\Lambda_{n+1}$ is the list $\Lambda_n$ plus any propositions of $M_{n+1}\setminus\{\emptyset,\Omega_{n+1}\}$ which are not/\emph{no more} listed in $\Lambda_n$\,.}
\end{itemize}
The first steps of the model construction are illustrated by a simple example in appendix~\ref{BayesMod:Comp&Ex}.
Appendix~\ref{BayesMod:Comp&Ex} also compares this construction with the Product Space CEA of~\cite{goodman2};
this comparison involves a discussion about the logical independence.
\subsubsection{Properties of $(\Omega_n, M_n, h_n, f_n,\Lambda_n)_{n\in\Nset}$}
\label{Omega:properties}
It is proved recursively:
\begin{description}
\item[$\rien\quad\bullet$] $\mu_n:M_n\rightarrow M_{n+1}$ is a one-to-one Boolean morphism,
\item[$\rien\quad\bullet$] If $A,B\in M_n$ and $f_n(B,A)$ is defined, then $f_{n+1}\bigl(\mu_n(B),\mu_n(A)\bigr)=\mu_n\bigl(f_n(B,A)\bigr)$\,,
\item[$\rien\quad\tilde\beta1$.] Let $A,B\in M_n$ such that $f_n(B,A)$ is defined.
\\Then $A\subset B$ and $A\ne\emptyset$ imply $f_n(B,A)=\Omega_n$\,,
\item[$\rien\quad\tilde\beta2$.] Let $A,B,C\in M_n$ such that $f_n(B,A)$, $f_n(C,A)$ and $f_n(B\cup C,A)$ are defined.
\\Then $f_n(B\cup C,A)= f_n(B,A)\cup f_n(C,A)$\,,
\item[$\rien\quad\tilde\beta3$.] Let $A,B\in M_n$ such that $f_n(B,A)$ is defined.
\\Then $A\cap f_n(B,A) = A\cap B$\,,
\item[$\rien\quad\tilde\beta4$.] Let $A,B\in M_n$ such that $f_n(B,A)$ and $f_n(\sim B,A)$ are defined.
\\Then $f_n(\sim B,A)=\sim f_n(B,A)$\,,
\item[$\rien\quad\tilde\beta6w$.] Let $A,B\in M_n$ such that $f_n(B,A)$ and $f_n(B,\sim A)$ are defined.
\\Then $f_n(B,A)=B$ implies $f_n(B,\sim A)=B$\,.
\end{description}
\emph{Proofs are given in appendix~\ref{Appendix:MainProof}}
\subsubsection{Limit}
The limit $\bigl(\Omega_\infty,M_\infty,h_\infty,f_\infty\bigr)$ is defined as follows:
\begin{itemize}
\item Set $\displaystyle\Omega_\infty=\left\{\left.(\omega_n|_{n\in\Nset})\in\prod_{n\in\Nset}\Omega_n\;\right/\;\forall n\in\Nset,\,\omega_{n+1}\in\mu_n(\omega_n)\right\}$\,;\\[10pt]
\emph{Useful definitions:}
\begin{itemize}
\item For any $n\in\Nset$ and any $A\in M_n$\,, $\displaystyle A_\infty=\bigl\{(\omega_k|_{k\in\Nset})\in\Omega_\infty
\;\big/\;\omega_n\in A\bigr\}$\,.
\emph{
The subset $A_\infty$ is a mapping of $A$ within $\Omega_\infty$\,.
It is noticed that this mapping is invariant with $\mu$, \emph{i.e.} $(A_{[m]})_\infty=A_\infty$ for $m>n$\,,
}
\item For any $n\in\Nset$\,, $\displaystyle M_{n:\infty}=\bigl\{A_\infty\;\big/\;A\in M_n\bigr\}$\,.
\emph{
The structure $M_{n:\infty}$ is an isomorphic mapping of the structure $M_n$ within $\Omega_\infty$\,.
It is noticed that $(M_{n:\infty}|_{n\in\Nset})$ is a monotonic sequence, \emph{i.e.} $M_{n:\infty}\subset M_{n+1:\infty}$\,,
}
\end{itemize}
\item Set $\displaystyle M_\infty=\bigcup_{n\in\Nset}M_{n:\infty}$\,,
\item Set $h_\infty(\theta)=\bigl(h_0(\theta)\bigr)_\infty$ for any $\theta\in\Theta$\,,
\item Let $A,B\in M_\infty$\,.
Then there is $n$ and $a,b\in M_n$ such that $A=a_\infty$\,, $B=b_\infty$ and $f_n(b,a)$ is defined (subsequent proposition). 
Set $f_\infty(B,A)=\bigl(f_{n}(b,a)\bigr)_\infty$\,.
\end{itemize}
This definition is justified by the following propositions:
\paragraph{Proposition 1.}
For any $n\in\Nset$, $m>n$ and $A\in M_n$, $(A_{[m]})_\infty=A_\infty$\,.
\begin{description}
\item[Proof.] By definition of $\Omega_\infty$ and since the $\mu_k|_{k\in\Nset}$ are one-to-one Boolean morphism:
$$
(A_{[m]})_\infty=\bigl\{(\omega_k|_{k\in\Nset})\in\Omega_\infty
\;\big/\;\omega_m\in A_{[m]}\bigr\}=\bigl\{(\omega_k|_{k\in\Nset})\in\Omega_\infty
\;\big/\;\omega_n\in A\bigr\}
=A_\infty\;.
$$
\item[$\Box\Box\Box$]\rien
\end{description}
\emph{Corollary} $M_{n:\infty}\subset M_{n+1:\infty}$\,.
\paragraph{Proposition 2.} $M_{n:\infty}$ is a Boolean subalgebra of $\mathcal{P}(\Omega_\infty)$ and is isomorph to $M_n$ by the morphism $A\mapsto A_\infty$.
\\[5pt]
\emph{Proof is obvious from the definition of $M_{n:\infty}$}\,.\\[5pt]
From now on, $M_n$ will be considered as a subalgebra of $\mathcal{P}(\Omega_\infty)$.
\paragraph{Proposition 3.} Let $A,B\in M_\infty$\,.
Then there is $n$ and $a,b\in M_n$ such that $A=a_\infty$\,, $B=b_\infty$ and $f_n(b,a)$ is defined.
\begin{description}
\item[Proof.]
Since $(M_{k:\infty}|_{k\in\Nset})$ is a monotonic sequence\,, there is $m\in\Nset$ such that $A,B\in M_{m:\infty}$\,.\\
Let $a,b\in M_m$ be such that $A=a_\infty$ and $B=b_\infty$\,.\\
By definition of the list $\Lambda_m$, there is $p\in[\![s_{m},f_{m}-1]\!]$ such that $\lambda_p=a$\,.\\
As a consequence, $f_{p+1}(b_{[p+1]},a_{[p+1]})$ exists.\\
But then hold $A=(a_{[p+1]})_\infty$ and $B=(b_{[p+1]})_\infty$\,.\\
Finally $n=p+1$ answers to the proposition.
\item[$\Box\Box\Box$]\rien
\end{description}
\paragraph{Proposition 4.} The definition of $f_\infty$ does not depend on the choice of $n$.
\begin{description}
\item[Proof.]
Let $n$, $m>n$, $a,b\in M_n$ and $c,d\in M_m$ such that $A=a_\infty=c_\infty$ and $B=b_\infty=d_\infty$\,.\\
Assume also that $f_n(b,a)$ and $f_m(d,c)$ exist.\\[3pt]
Then $(a_{[m]})_{\infty}=c_\infty$ and $(b_{[m]})_{\infty}=d_\infty$.\\
Since $M_n$ and $M_{n:\infty}$ are isomorph, it follows $a_{[m]}=c$ and $b_{[m]}=d$\,.\\
But it has been shown in a previous section that $f_m(b_{[m]},a_{[m]})=f_n(b,a)_{[m]}$\,.\\
Finally $\bigl(f_n(b,a)\bigr)_\infty=\bigl(f_m(d,c)\bigr)_\infty$\,.
\item[$\Box\Box\Box$]\rien
\end{description}
\paragraph{Proposition 5.} $\bigl(\Omega_\infty,M_\infty,h_\infty,f_\infty\bigr)$  verifies $\beta1$, $\beta2$, $\beta3$, $\beta4$ and $\beta6w$.\\[5pt]
The properties $\beta\ast$ are inherited from $M_n,f_n|_{n\in\Nset}$\,, by means of the properties $\tilde\beta\ast$.
\paragraph{Conclusion.} \underline{$\bigl(\Omega_\infty,M_\infty,h_\infty,f_\infty\bigr)$ is a conditional model for DmBL$_\ast$.}
\subsubsection{Implied Kripke model for DmBL$_\ast$}
By means of the \emph{Model transfer} property, a Kripke model for DmBL$_\ast$ is derived from $\bigl(\Omega_\infty,M_\infty,h_\infty,f_\infty\bigr)$.
\underline{This Kripke model is denoted $\mathcal{B}=\bigl(\Omega_\infty,R_{\mathcal{B}},H_{\mathcal{B}},f_\infty\bigr)$\,.}
\subsubsection{Partial Completeness}
It is above the scope of this work to construct a complete model for DmBL$_\ast$.
However, it is shown here that $\mathcal{B}$ is sufficiently complete to characterize the operator $(|)$ in DmBL$_\ast$.
\paragraph{Proposition 1.}
By construction, $(\Omega_\infty,H_{\mathcal{B}})$ is a complete model for the classical logic $C$, when $H_{\mathcal{B}}$ is restricted to the propositions of $\mathcal{L}_C$.
\paragraph{Proposition 2.}
Let $\phi\in\mathcal{L}$ be a proposition constructed \underline{without $\Box$ or $\Diamond$}\,.
Then $\vdash\phi$ in DmBL$_\ast$ if and only if $H_{\mathcal{B}}(\phi)=\Omega_\infty$\,.
\\[5pt]
Proof is done in appendix~\ref{Appendix:ProofOfAlmostCompletude}\,.
\subsection{Coherence properties}
The model $\mathcal{B}$ clearly shows that DmBL$_\ast$ is coherent.
It also demonstrates that the conditional operator $(|)$ is not trivial.
Since $(\Omega_\infty,H_{\mathcal{B}})$ is a complete model for $C$, DmBL$_\ast$ is also an extension of the classical logic:
$\vdash\phi$ implies $\vdash_C\phi$.
But a stronger property holds:
\paragraph{Non distortion.}
Let $\phi$ be a classical proposition.
Assume that $\vdash\Box\phi\vee\Box\neg\phi$ in DmBL$_\ast$.
Then $\vdash_C\phi$ or $\vdash_C\neg\phi$.
\\[5pt]
\emph{Interpretation: DmBL$_\ast$ does not ``distort'' the \emph{classical} propositions.
More precisely, a property like $\vdash\Box\phi\vee\Box\neg\phi$ would add some knowledge about $\phi$, since it says that either $\phi$ or $\neg\phi$ is ``sure''.
But the \emph{non distortion} just tells that such property is impossible unless there is a trivial knowledge about $\phi$ within the classical logic.}
\begin{description}
\item[Proof.]
Assume $\vdash\Box\phi\vee\Box\neg\phi$\,.\\
Since $\mathcal{B}$ is a model for DmBL$_\ast$, it comes $H_{\mathcal{B}}(\Box\phi\vee\Box\neg\phi)=\Omega_\infty$\,.\\
Then $H_{\mathcal{B}}(\Box\phi)\cup H_{\mathcal{B}}(\Box\neg\phi)=\Omega_\infty$\,.\\
Since $H_{\mathcal{B}}(\Box \psi)=\emptyset$ or $\Omega_\infty$ for any $\psi\in\mathcal{L}$\,, it comes $H_{\mathcal{B}}(\Box\phi)=\Omega_\infty$ or $H_{\mathcal{B}}(\Box\neg\phi)=\Omega_\infty$\,.\\
At last, $H_{\mathcal{B}}(\phi)=\Omega_\infty$ or $H_{\mathcal{B}}(\neg\phi)=\Omega_\infty$\,.\\
But $(\Omega_\infty,H_{\mathcal{B}})$ is a complete Boolean model for $C$\,, which implies $\vdash_C\phi$ or $\vdash_C\neg\phi$\,.
\item[$\Box\Box\Box$]\rien
\end{description}
Another \emph{non distortion} property will be derived subsequently, by the probabilistic extension.
\section{Probabilistic extension}
\label{Section:Proba:DmBL}
\subsection{Probability over propositions,} a minimal$^\dagger$ definition.\\
{\small $\dagger$ This definition is related to finite probabilities and excludes any Bayesian consideration.}
\\[7pt]
Probabilities are classically defined over measurable sets.
However, this is only a manner to model the notion of probability, which is essentially an additive measure of the belief of logical propositions \cite{paass}.
Probability could be defined without reference to the measure theory, at least when the propositions are countable.
The notion of probability is explained now within a strict propositional formalism.
Conditional probabilities are excluded from this definition, but the notion of independence is considered.
\vspace{5pt}\\
Intuitively, a probability over a set of logical propositions is a measure of belief which is additive (disjoint propositions are adding their chances) and increasing with the propositions.
This measure should be zeroed for the ever-false propositions and full for the ever-true propositions.
Moreover, \emph{a probability is a multiplicative measure for independent propositions}.
These intuitions are now formalized.
\paragraph{Definition for classical propositions.}
A probability $\pi$ over $C$ is a $\Rset^+$ valued function such that for any proposition $\phi$ and $\psi$ of $\mathcal{L}_C$\,:
\begin{description}
\item[\rien$\quad$\emph{Equivalence.}]$\phi\equiv_C\psi$ implies $\pi(\phi)=\pi(\psi)$\,,
\item[\rien$\quad$\emph{Additivity.}]$\pi(\phi\wedge\psi)+\pi(\phi\vee\psi)=\pi(\phi)+\pi(\psi)$\,,
\item[\rien$\quad$\emph{Coherence.}]$\pi(\bot)=0$\,,
\item[\rien$\quad$\emph{Finiteness.}]$\pi(\top)=1$\,.
\end{description}
\subparagraph{Property.}
The coherence and additivity implies the increase of $\pi$:
\begin{description}
\item[\rien$\quad$\emph{Increase.}]$\pi(\phi\wedge\psi)\le \pi(\phi)$\,.
\end{description}
\begin{description}
\item[Proof.] 
Since $\phi\equiv_C(\phi\wedge\psi)\vee(\phi\wedge\neg\psi)$ and $(\phi\wedge\psi)\wedge(\phi\wedge\neg\psi)\equiv_C\bot$, the additivity implies:
$$
\pi(\phi)+\pi(\bot)=\pi(\phi\wedge\psi)+\pi(\phi\wedge\neg\psi)\;.
$$
From the coherence $\pi(\bot)=0$\,,
it is deduced $\pi(\phi)=\pi(\phi\wedge\psi)+\pi(\phi\wedge\neg\psi)$\,.\\
Since $\pi$ is non negatively valued, $\pi(\phi)\ge \pi(\phi\wedge\psi)$\,.
\item[$\Box\Box\Box$]\rien
\end{description}
\paragraph{Definition for DmBL/DmBL$_\ast$.}
In this case, we have to deal with independence notions.\\[5pt]
A probability $P$ over DmBL/DmBL$_\ast$ is a $\Rset^+$ valued function, which verifies (replace $\equiv_C$ by $\equiv$ and $\pi$ by $P$) \emph{equivalence}, \emph{additivity}, \emph{coherence}, \emph{finiteness} and:
\begin{description}
\item[\rien$\quad$\emph{Multiplicativity.}]$\vdash\phi\times\psi$ implies $P(\phi\wedge\psi)=P(\phi)P(\psi)$\,.
\end{description}
for any proposition $\phi$ and $\psi$ of $\mathcal{L}$\,. 
\subsection{Probability extension over DmBL$_\ast$}
\label{ProbExt:wDmBL}
\paragraph{Property.}
Let $\pi$ be a probability defined over $C$\,, the classical logic, such that $\pi(\phi)>0$ for any $\phi\not\equiv_C\bot$.
Then, there is a (multiplicative) probability $P$ defined over DmBL$_\ast$ such that $P(\phi)=\pi(\phi)$ for any classical proposition $\phi\in\mathcal{L}_C$\,.
\\\\
\emph{Remark: this is another non distortion property, since the construction of DmBL$_\ast$ puts no constraint over probabilistic classical propositions.}
\\[5pt]
Proof is done in appendix~\ref{Appendix:Probabilition}.
\paragraph{Corollary.}
Let $\pi$ be a probability defined over $C$\,.
Then, there is a probability $P$ defined over DmBL$_\ast$ such that $P(\phi)=\pi(\phi)$ for any $\phi\in\mathcal{L}_C$\,.
\begin{description}
\item[Proof.]
Let $\Sigma=\left\{\left.\bigwedge_{\theta\in\Theta}\epsilon_\theta\;\right/\;\epsilon\in\prod_{\theta\in\Theta}\{\theta,\neg\theta\}\right\}$\,.\\
For any real number $e >0$\,, define then $\pi_e $ the probability over $\mathcal{L}_C$ such that:
$$
\forall\sigma\in\Sigma\,,\; \pi_e (\sigma)=\frac{e }{\mathrm{card}(\Sigma)}+(1-e )\pi(\sigma)
\;. 
$$
Let $P_e $ be the probability over DmBL$_\ast$ constructed from $\pi_e $ as defined in appendix~\ref{Appendix:Probabilition}.\\
By~\ref{AppC:Conclude}\,, there is a rational function $R_\phi$ such that $P_e (\phi)=R_\phi(e )$ for any $\phi\in\mathcal{L}$\,.\\
Now $0\le R_\phi(e )\le 1$\,;
since $R_\phi(e )$ is rational and bounded, $\lim_{e \rightarrow 0+}R_\phi(e )$ exists.\\
Define $P(\phi)=\lim_{e \rightarrow 0+}R_\phi(e )$\,, for any $\phi\in\mathcal{L}$.\\
The additivity, coherence, finiteness and multiplicativity are obviously inherited by $P$.\\
At last, it is clear that $P(\sigma)=\pi(\sigma)$ for any $\sigma\in\Sigma$\,.
\item[$\Box\Box\Box$]\rien
\end{description}
\subsection{Model and probability extension for DmBL}
Let $\mathcal{K}$ be the set of all (multiplicative) probabilities $P$ over DmBL$_\ast$ such that $P(\phi)>0$ for any $\phi\not\equiv\bot$\,,
and define the sequences $\mathcal{K}(\phi)=(P(\phi))_{P\in\mathcal{K}}$ for any $\phi\in\mathcal{L}$\,.\\
Then set $
\mathcal{L}_{\mathcal{K}}=\mathcal{K}(\mathcal{L})=\bigl\{\mathcal{K}(\phi)\;\big/\;\phi\in\mathcal{L}\bigr\}\;,
$ a subset of ${\Rset^+}^{\mathcal{K}}$\,.\\
The operators $\neg$, $\wedge$ and $(|)$ are canonically implied over $\mathcal{L}_{\mathcal{K}}$\,:\footnote{The operators $\vee$ and $\rightarrow$ are derived from $\wedge$ and $\neg$ as usually; modalities are not considered.}
$$
\neg\mathcal{K}(\phi)=\mathcal{K}(\neg\phi)\,,\ 
\mathcal{K}(\phi)\wedge\mathcal{K}(\psi)=\mathcal{K}(\phi\wedge\psi)
\quad\mbox{and}\quad
\bigl(\mathcal{K}(\psi)\big|\mathcal{K}(\phi)\bigr)=
\mathcal{K}\bigl((\psi|\phi)\bigr)\;.
$$
Since any $P\in\mathcal{K}$ verifies the equivalence property, it comes $\mathcal{K}(\phi)=\mathcal{K}(\psi)$ when $\phi\equiv\psi$ in DmBL$_\ast$.
As a direct consequence, $\bigl(\mathcal{L}_{\mathcal{K}},\neg,\wedge,(|)\bigr)$ is a conditional-like model of DmBL$_\ast$ (the structure is a Boolean algebra but not derived from set operators. that is the only difference with conditional models).
\paragraph{Property.}
$\bigl(\mathcal{L}_{\mathcal{K}},\neg,\wedge,(|)\bigr)$ is a conditional-like model of DmBL.
\begin{description}
\item[Proof.]Let $P\in\mathcal{K}$\,; $P$ is multiplicative.\\
Since $\vdash(\psi|\phi)\times\phi$ and $(\psi|\phi)\wedge\phi\equiv\psi\wedge\phi$ in DmBL$_\ast$, it comes $P\bigl((\psi|\phi)\bigr)P(\phi)=P(\psi\wedge\phi)$\,.
\\[3pt]
Now assume $\bigl(\mathcal{K}(\psi)\big|\mathcal{K}(\phi)\bigr)=\mathcal{K}(\psi)$\,, with $\psi\not\equiv\bot$\,.\\
Then $\mathcal{K}\bigl((\psi|\phi)\bigr)=\mathcal{K}(\psi)$, and $P\bigl((\psi|\phi)\bigr)=P(\psi)$ for any $P\in\mathcal{K}$\,.\\
Then $P(\phi)=\frac{P(\psi\wedge\phi)}{P\bigl((\psi|\phi)\bigr)}=\frac{P(\psi\wedge\phi)}{P(\psi)}=P\bigl((\phi|\psi)\bigr)$ for any $P\in\mathcal{K}$\,,
\\
and $\bigl(\mathcal{K}(\phi)\big|\mathcal{K}(\psi)\bigr)=\mathcal{K}\bigl((\phi|\psi)\bigr)=\mathcal{K}(\phi)$\,.\\
Since moreover $(\phi|\bot)\equiv\phi$ and $(\bot|\phi)\equiv\bot$ in DmBL$_\ast$, the model verifies $\beta6$\,.
\item[$\Box\Box\Box$]\rien
\end{description}
Notice that it was only needed the equivalence and multiplicative properties for the elements of $\mathcal{K}$\,.
It is thus possible to construct a more general model by relaxing $\mathcal{K}$\,.
\paragraph{Probability extension.}
For any $\mathcal{K}(\phi)\in\mathcal{L}_{\mathcal{K}}$ and any $P\in\mathcal{K}$\,, define the $\Rset^+$-valued mapping ${\hat P}\bigl(\mathcal{K}(\phi)\bigr)=P(\phi)$ (this mapping, a projection, is indeed well defined).\\
By construction, ${\hat P}$ is naturally a multiplicative probability over $\mathcal{L}_{\mathcal{K}}$\,.
Moreover, the probabilistic extensions defined in appendix~\ref{Appendix:Probabilition} are also elements of $\mathcal{K}$\,.
As a consequence, the deductions of section~\ref{ProbExt:wDmBL} are still working for $\mathcal{L}_{\mathcal{K}}$\,.
The extension property is thus derived:
\begin{quote}
Let $\pi$ be a probability defined over $C$\,.
Then, there is a probability $P$ defined over DmBL such that $P(\phi)=\pi(\phi)$ for any $\phi\in\mathcal{L}_C$\,.
\end{quote}
\paragraph{Non distortion.}
Let $\phi$ be a classical proposition.
Assume that $\vdash\Box\phi\vee\Box\neg\phi$ in DmBL.
Then $\vdash_C\phi$ or $\vdash_C\neg\phi$.
\begin{description}
\item[Proof.]
Consider the Kripke model for DmBL derived from the conditional model $\bigl(\mathcal{L}_{\mathcal{K}},\neg,\wedge,(|)\bigr)$\,.
\\
In this model, the value of $H(\Box\phi)$ is either $\mathcal{K}(\bot)$ or $\mathcal{K}(\top)$\,.\\
Then, $H(\Box\phi\vee\Box\neg\phi)=\mathcal{K}(\top)$ implies $H(\Box\phi)=\mathcal{K}(\top)$ or $H(\Box\neg\phi)=\mathcal{K}(\top)$\,.\\
Then $H(\phi)=\mathcal{K}(\top)$ or $H(\neg\phi)=\mathcal{K}(\top)$\,.\\
It follows $\forall P\in\mathcal{K}\,,\; P(\phi)=1$ or $\forall P\in\mathcal{K}\,,\; P(\neg\phi)=1$\,, and by the probabilistic extension:
$\forall \pi\,,\; \pi(\phi)=1$ or $\forall \pi\,,\; \pi(\neg\phi)=1$\,, where $\pi$ is any probability over $C$\,.\\
At last, $\vdash_C\phi$ or $\vdash_C\neg\phi$\,.
\item[$\Box\Box\Box$]\rien
\end{description}
\subsection{Properties of the conditional}
\paragraph{Bayes inference.}
Assume a probability $P$ defined over DmBL/DmBL$_\ast$.
Define $P(\psi|\phi)$ as an abbreviation for $P\bigl((\psi|\phi)\bigr)$\,.
Then:
$$
P(\psi|\phi)P(\phi)=P(\phi\wedge\psi)\;.
$$
\begin{description}
\item[Proof.]A consequence of \mbox{$(\psi|\phi)\wedge\phi\equiv\phi\wedge\psi$} and \mbox{$\vdash(\psi|\phi)\times\phi$}\,.
\item[$\Box\Box\Box$]\rien
\end{description}
As a corollary, it is also deduced
\underline{$P\bigl((\psi|\phi)\big|(\eta|\zeta)\bigr)P(\eta|\zeta)=P\bigl((\psi|\phi)\wedge(\eta|\zeta)\bigr)$}\,.
It is recalled that the closure of DGNW fails on this relation.
\paragraph{About Lewis' negative result.}
\subparagraph{The theorem.}
Lewis triviality may be expressed as follows:
\begin{quote}
Let $M_0$ be the set of measurable subsets of $\Omega_0$, and let $A,B\in M_0$ with $\emptyset\subsetneq B\subsetneq 
A\subsetneq \Omega_0$\,.
Then, it is impossible to construct a proposition $(B|A)\in M_0$ such that $\pi\bigl((B|A)\bigr)=\pi(B|A)\stackrel{\Delta}{=}\frac{\pi(A\cap B)}{\pi(A)}$ for any $\pi\in\mathcal{P}(M_0)$ with $0<\pi(B)<\pi(A)<1$\,
\end{quote}
\begin{description}
\item[Proof.]
For any propositions $C,D$\,, define $\pi_C(D)=\pi(D|C)=\frac{\pi(C\cap D)}{\pi(C)}$\,, when $\pi(C)>0$.\\
The proof of Lewis's result relies of the following calculus:
\begin{equation}
\label{Eq:Lewis}
\pi((E|D)|C)=\frac{\pi_C(D\cap E)}{\pi_C(D)}=\frac{\frac{\pi(C\cap D\cap E)}{\pi(C)}}{\frac{\pi(D\cap C)}{\pi(C)}}= \frac{\pi(C\cap D\cap E)}{\pi(D\cap C)}=\pi(E|C\cap D)\;.
\end{equation}
Then, it is inferred (with the notation $\sim B=\Omega_0\setminus B$):
$$
\begin{array}{@{}l@{}}\displaystyle
\frac{\pi(A\cap B)}{\pi(A)}=\pi(B|A)=\pi((B|A)|B)\pi(B)+\pi((B|A)|\sim B)\pi(\sim B)
\vspace{3pt}\\\displaystyle
\hspace{30pt}=\pi(B|B\cap A)\pi(B)+\pi(B|\sim B\cap A)\pi(\sim B)=1\times \pi(B)+0\times \pi(\sim B)=\pi(B)\;.
\end{array}
$$
As a consequence, $A$ and $B$ are probabilistically independent for any possible choice of $\pi$ such that $0<\pi(B)<\pi(A)<1$\,.
This is impossible.
\item[$\Box\Box\Box$]\rien
\end{description}
\subparagraph{Discussion.}
What is refuted by Lewis is the hypothesis $(B|A)\in M_0$, with $\emptyset\subsetneq B\subsetneq 
A\subsetneq \Omega_0$\,.
On the contrary, the hypotheses $(\emptyset|A)\in M_0$\,, $(\Omega_0|A)\in M_0$\,, $(B|\Omega_0)\in M_0$\,, $(B|\emptyset)\in M_0$ and $(C|A)\in M_0$ for $A\subset C$, are not refuted by Lewis result itself, and may be used in the construction of conditional.
This is an argument for removing trivial conditionings.
\\[3pt]
Now, assume a model of DmBL to be constructed, and let $M_0$ be the Boolean algebra associated to the classical propositions (\emph{i.e.} the unconditioned propositions).
Our model should be compatible with any probability defined over the unconditioned propositions (this is a minimal requirement!)
Lewis's triviality then make necessary the following hypothesis in the model construction:
\emph{$(B|A)\not\in M_0$ for any $A,B\in M_0$ such that $\emptyset\subsetneq B\subsetneq 
A\subsetneq \Omega_0$\,.}
Since DmBL implies $(\psi|\phi)\wedge\phi\equiv\phi\wedge\psi$\,,
it is inferred $(\psi|\phi)\equiv (\phi\wedge\psi)\vee\bigl(\neg\phi\wedge(\psi|\phi)\bigr)$ and:
\begin{center}
\emph{$(B|A)\in M_0$ if and only if $(\Omega_0\setminus A)\cap(B|A)\in M_0$\,, for any $A,B\in M_0$\,.}
\end{center}
At last, a model of DmBL have to conform to: \underline{$(\Omega_0\setminus A)\cap(B|A)\not\in M_0$\,, for any $A,B\in M_0$} \underline{such that $\emptyset\subsetneq B\subsetneq 
A\subsetneq \Omega_0$\,.}
A thorough look at case 1 of the model construction in section~\ref{Section:Model:DmBL}, clearly shows that our model is compliant with such condition.
But it is noteworthy that the Product Space CEA is also compliant with it (refer to \cite{goodman2} or to appendix \ref{BayesMod:Comp&Ex}).
This condition is probably one of the keys for constructing a model for Bayesian.
\subparagraph{And DmBL\dots}
The previous extension theorems has shown that for any probability $\pi$ defined over $C$\,, it is possible to construct a probability $P$ over DmBL which extends $\pi$\,.
This result by itself shows that DmBL avoids Lewis triviality.
The question is why?
In fact the coexistence of Lewis and the extension theorem infers necessarily the following conditions:
\begin{itemize}
\item Define $Q_\phi$ by $Q_\phi(\psi)=P\bigl((\psi|\phi)\bigr)$\,.
Then $Q_\phi$ is \emph{additive}, \emph{coherent} and \emph{finite}\,, which makes it a probability in the classical meaning.
But \underline{$Q_\phi$ is not necessarily \emph{multiplicative}}\,.
\item Let $\phi,\psi\in\mathcal{L}_C$ be such that $\vdash\psi\rightarrow\phi$\,, $\psi\not\equiv\bot$\,, $\psi\not\equiv\phi$ and $\phi\not\equiv\top$\,.
\\
Then $\forall\eta\in\mathcal{L}_C\,,\neg\phi\wedge(\psi|\phi)\not\equiv\eta$\,.
\end{itemize}
Conversely, both conditions are sufficient for the coexistence of Lewis and the extension.
\\[5pt]
The last condition is the necessary construction condition, implied by the discussion.
\\[3pt]
The first condition expresses a limitation of DmBL.
This limitation is unavoidable: otherwise the derivation~(\ref{Eq:Lewis}) is possible, even if $(\psi|\phi)$ is not equivalent to a classical proposition.
\\[3pt]
Besides, assume $\phi\in\mathcal{L}_C$ and define the probability $\pi_\phi$ over $C$ by $\pi_\phi(\psi)=\pi(\phi\wedge\psi)/\pi(\phi)$\,.
Let $P_\phi$ be the extension of $\pi_\phi$ over DmBL.
Then $P_\phi\ne Q_\phi$\,.
\paragraph{Logical independence and probabilistic independence.}
The logical independence is a property \emph{stronger} than the probabilistic independence.
The independence $\vdash\psi\times\phi$ is equivalent to $(\psi|\phi)\equiv\psi$\,.
It is independent of the choice of a probability.
In the probabilist paradigm, $\phi$ and $\psi$ are said to be independent when $P(\psi|\phi)=P(\psi)$\,.
It is dependent of $P$.
In particular, the condition $\sigma\equiv\psi$ implies $P(\sigma)=P(\psi)$ but the converse is false.
The independence is a special case, where $\sigma=(\psi|\phi)$\,.
\\[5pt]
The regularity property is a particular illustration of this difference:
\begin{center}
$\vdash(\phi\wedge\eta)\rightarrow(\psi\wedge\eta)$\,, $\vdash\phi\times\eta$, $\vdash\psi\times\eta$ and $\vdash\Diamond \eta$
imply $\vdash\phi\rightarrow\psi$\,.
\end{center}
Now, considering $\phi,\psi,\eta$ as measurable sets, the probabilistic independences $P(\phi\cap\eta)=P(\phi)P(\eta)$\,, $P(\psi\cap\eta)=P(\psi)P(\eta)$ and the inclusion $\phi\cap\eta\subset\psi\cap\eta$ do not yield the inclusion $\phi\subset\psi$ in general, even when $\eta\ne\emptyset$.
\section{Conclusion}
\label{Fus2004::Sec:8}
In this contribution, a conditional logic, DmBL, has been defined and studied.
This Logic is essentially different to the VCU of Lewis, since it implies a classical nature to the sub-universe.
By the way, it shares common properties with notable Conditional Event Algebras.
This logic has been proved to be coherent.
A complete model has been defined for the propositions and the conditionals (not the modalities), in the case of the weakened logic DmBL$_\ast$.
It has been shown that any probability over the classical propositions could be extended to DmBL/DmBL$_\ast$, in compliance with the independence relation.
Then, the probabilistic Bayesian rule has been recovered from DmBL/DmBL$_\ast$.
\\[5pt]
There are still many open questions.
For example, it is certainly possible to bring some enrichment to the conditional of DmBL, by means of additional axioms.
Is it possible to recover some specific equivalences of DGNW, or other CEA?
\\[5pt]
From a strictly logical viewpoint, the Deterministic modal Bayesian Logic has also some interesting properties.
For example, the notion of independence in DmBL have some surprising logical consequences in the deductions (\emph{e.g.} regularity with an inference).
Is there a possible application in the constructions of theories?
\\[5pt]
At last it is recalled our first motivation, and this work will be applied in a next future to the logical comparison of the Bayesian with other theories of uncertain informations.
\small
%

%
\appendix
\section{Proof: properties of $(\Omega_n, M_n, h_n, f_n,\Lambda_n)_{n\in\Nset}$}
\label{Appendix:MainProof}
To be proved:
\emph{\begin{description}
\item[$\rien\quad\bullet_\mu$] $\mu_n:M_n\rightarrow M_{n+1}$ is a one-to-one Boolean morphism\,,
\item[$\rien\quad\bullet_f$] $f_{n+1}\bigl(\mu_n(B),\mu_n(A)\bigr)=\mu_n\bigl(f_n(B,A)\bigr)$\,,
\item[$\rien\quad\tilde\beta1$.] $A\subset B$ and $A\ne\emptyset$ imply $f_n(B,A)=\Omega_n$\,,
\item[$\rien\quad\tilde\beta2$.] $f_n(B\cup C,A)= f_n(B,A)\cup f_n(C,A)$\,,
\item[$\rien\quad\tilde\beta3$.] $A\cap f_n(B,A) = A\cap B$\,,
\item[$\rien\quad\tilde\beta4$.] $f_n(\sim B,A)=\sim f_n(B,A)$\,,
\item[$\rien\quad\tilde\beta6w$.] $f_n(B,A)=B$ implies $f_n(B,\sim A)=B$\,,
\end{description}
being assumed $A,B,C\in M_n$\,, and $f_n(\cdot,\cdot)$ defined for the considered cases.
}
\\[5pt]
The proof is recursive and needs to consider the three cases in the definition of $(\mu_n,f_n)$\,.
\\[3pt]
The properties $\tilde\beta\ast$ are obvious for $n=0$, since $f_0$ is only defined by $f_(A,\emptyset)=f_0(A,\Omega_0)=A$\,.
From now on, it is assumed that $\tilde\beta\ast$ hold true for $k\le n$, and that $\bullet_\mu$ and $\bullet_f$ hold true for $k\le n-1$\,.
The subsequent paragraphs establish the proof of $\beta\ast$ for $n+1$ and the proof of $\bullet_\mu$ and $\bullet_f$ for $n$.
\paragraph{Preliminary remark.} It is noticed that $\tilde\beta2$ and $\tilde\beta4$ imply:
$$\tilde\beta7:\ f_n(B\cap C,A)= f_n(B,A)\cap f_n(C,A)\,.$$
\subsection{Lemma.}
\label{lemma:include}
$\bigcup_{i\in I_n}\Pi_n(i)=b_n$
and
$\bigcup_{i\in I_n}\Gamma_n(i)=\sim b_n$\,;
in particular, $\Pi_n(i)\cap\Gamma_n(j)=\emptyset$ for any $i,j\in I_n$\,.\\[5pt]
Moreover $\Pi_n(i)\cap\Pi_n(j)=\Gamma_n(i)\cap\Gamma_n(j)=\emptyset$ for any $i,j\in I_n$ such that $i\ne j$\,.
\begin{description}
\item[Proof.]The proof is obvious for case 1.\\
Now, let consider case 0.\\
By definition $
\bigcup_{i\in I_n}\Pi_n(i)=\left(\bigcup_{\omega\in\mu_\nu(b_\nu)}\omega_{[n]}\right)
\cap
\left(\bigcup_{\omega'\in\sim\mu_\nu(b_\nu)}f_n(\omega'_{[n]},\sim b_n)\right)
$\,.\\
By recursion hypothesis over $\tilde\beta1$ it comes $f_n(\sim b_n,\sim b_n)=\Omega_n$\,.
\\
Then by $\tilde\beta2$\,, $\bigcup_{i\in I_n}\Pi_n(i)=b_n\cap\Omega_n= b_n\;.$\\[5pt]
For any $\omega_1,\omega_2\in \mu_\nu(b_\nu)$ such that $\omega_1\ne\omega_2$\,, it comes $\tilde\beta7$ by $\tilde\beta2,\tilde\beta4$, and:
$$
f_n(\omega_{1[n]},b_n)\cap f_n(\omega_{2[n]},b_n)=f_n(\omega_{1[n]}\cap\omega_{2[n]},b_n)=f_n(\emptyset,b_n)=\emptyset\;.
$$
Finally $\Gamma_n(i)\cap\Gamma_n(j)=\emptyset$ for any $i,j\in I_n$ such that $i\ne j$\,.\\[5pt]
The results are similarly proved for $\Gamma_n$\,.
\item[$\Box\Box\Box$]\rien
\end{description}
\emph{Corollary 1.}
$$
\mu_n(b_n)=T(\sim\mu_n(b_n))=\bigcup_{i\in I_n}\Pi_n(i)\times\Gamma_n(i)
\quad\mbox{and}\quad
\sim\mu_n(b_n)=T(\mu_n(b_n))=\bigcup_{i\in I_n}\Gamma_n(i)\times\Pi_n(i)\;.
$$
\emph{Corollary 2.}
$$f_{n+1}\bigl(C,\mu_n(b_n)\bigr)=\bigl(C\cap\mu_n(b_n)\bigr)\cup\bigl(T(C)\cap\sim\mu_n(b_n)\bigr)$$
and
$$f_{n+1}\bigl(C,\sim\mu_n(b_n)\bigr)=\bigl(T(C)\cap\mu_n(b_n)\bigr)\cup\bigl(C\cap\sim\mu_n(b_n)\bigr)\;.$$
Both corollary are obvious from the definition.
\subsection{Proof of $\bullet_\mu$}
\label{Proof:Cap}
The following properties (whose proofs are immediate) will be useful:
\begin{description}
\item[$\ell1.$] $(A\cup B)\times C=(A\times C)\cup (B\times C)$
and $A\times (B\cup C)=(A\times B)\cup (A\times C)$\,, for any $A,B,C$\,,
\item[$\ell2.$] $(A\cap B)\times C=(A\times C)\cap (B\times C)$
and $A\times (B\cap C)=(A\times B)\cap (A\times C)$\,, for any $A,B,C$\,,
\item[$\ell3.$] $C\cap D=\emptyset$ implies $(C\times A)\cap(D\times B)=(A\times C)\cap(B\times D)=\emptyset$\,, for any $A,B,C,D$\,,
\item[$\ell4.$] $(A\cup B)\cap (C\cup D)=\emptyset$ implies $(A\cap B)\cup(C\cap D)=(A\cup C)\cap(B\cup D)$\,, for any $A,B,C,D$\,.
\item[$\ell5.$] $(A\cup B)\cap (C\cup D)=\emptyset$ and $A\cup C=B\cup D$ imply $A=B$ and $C=D$\,, for any $A,B,C,D$\,.
\item[$\ell6.$] $C\cap D=\emptyset$ implies $(C\cup D)\setminus\bigl((A\cap C)\cup(B\cap D)\bigr)=(C\setminus A)\cup(D\setminus B)$\,, for any $A,B,C,D$\,.
\end{description}
\paragraph{\underline{Proof of $\mu_n(\Omega_n)=\Omega_{n+1}$ and $\mu_n(\emptyset)=\emptyset$\,.}}
Immediate from the definitions.
\paragraph{\underline{Proof of $\mu_n(A\cap B)=\mu_n(A)\cap\mu_n(B)$\,.}}
By applying $\ell2$, it is deduced:
$$\begin{array}{@{}l@{}}
\mu_n(A\cap B)=\bigcup_{i\in I_n} \biggl(\Bigl(\bigl(A\cap\Pi_n(i)\bigr)\times\Gamma_n(i)\Bigr)
\cap
\Bigl(\bigl(B\cap\Pi_n(i)\bigr)\times\Gamma_n(i)\Bigr)\biggr)
\vspace{5pt}\\\rien\hspace{50pt}
\cup\ 
\bigcup_{i\in I_n}
\biggl(\Bigl(\bigl(A\cap\Gamma_n(i)\bigr)\times\Pi_n(i)\Bigr)
\cap
\Bigl(\bigl(B\cap\Gamma_n(i)\bigr)\times\Pi_n(i)\Bigr)\biggr)
\;.
\end{array}$$
By lemma~\ref{lemma:include}, and applying $\ell3$ and $\ell4$, it is deduced:
$$\begin{array}{@{}l@{}}
\mu_n(A\cap B)=\bigcup_{i\in I_n} \biggl(\Bigl(\bigl(A\cap\Pi_n(i)\bigr)\times\Gamma_n(i)\Bigr)
\cup
\biggl(\Bigl(\bigl(A\cap\Gamma_n(i)\bigr)\times\Pi_n(i)\Bigr)
\vspace{5pt}\\\rien\hspace{50pt}
\cap\ 
\bigcup_{i\in I_n}
\Bigl(\bigl(B\cap\Pi_n(i)\bigr)\times\Gamma_n(i)\Bigr)\biggr)
\cup
\Bigl(\bigl(B\cap\Gamma_n(i)\bigr)\times\Pi_n(i)\Bigr)\biggr)
=\mu_n(A)\cap\mu_n(B)\;.
\end{array}$$
\paragraph{\underline{Proof of $\mu_n(A\cup B)=\mu_n(A)\cup\mu_n(B)$\,.}}
Obviously deduced from $\ell1$.
\paragraph{\underline{$\mu_n$ is one-to-one.}}
Assume $\mu_n(A)=\mu_n(B)$\,;
then:
$$\begin{array}{@{}l@{}}
\bigcup_{i\in I_n} \biggl(\Bigl(\bigl(A\cap\Pi_n(i)\bigr)\times\Gamma_n(i)\Bigr)
\cup
\biggl(\Bigl(\bigl(A\cap\Gamma_n(i)\bigr)\times\Pi_n(i)\Bigr)
\vspace{5pt}\\\rien\hspace{50pt}
=
\bigcup_{i\in I_n}
\Bigl(\bigl(B\cap\Pi_n(i)\bigr)\times\Gamma_n(i)\Bigr)\biggr)
\cup
\Bigl(\bigl(B\cap\Gamma_n(i)\bigr)\times\Pi_n(i)\Bigr)\biggr)
\;.
\end{array}$$
By lemma~\ref{lemma:include}, and applying $\ell2$, $\ell3$ and $\ell5$, it is deduced for any $i\in I_n$\,:
$$
\bigl(A\cap\Pi_n(i)\bigr)\times\Gamma_n(i)
=
\bigl(B\cap\Pi_n(i)\bigr)\times\Gamma_n(i)
\ 
\mbox{ and }
\ 
\bigl(A\cap\Gamma_n(i)\bigr)\times\Pi_n(i)
=
\bigl(B\cap\Gamma_n(i)\bigr)\times\Pi_n(i)
\;.$$
Finally $A\cap\Pi_n(i)=B\cap\Pi_n(i)$ and $A\cap\Gamma_n(i)=B\cap\Gamma_n(i)$ for any $i\in I_n$\,, and:
$$
A\cap\bigcup_{i\in I_n}\bigl(\Pi_n(i)\cup\Gamma_n(i)\bigr)=B\cap\bigcup_{i\in I_n}\bigl(\Pi_n(i)\cup\Gamma_n(i)\bigr)\;.
$$
$A=B$ is deduced by applying the lemma.
\paragraph{Conclusion.} The previous results imply that $\mu_n$ is a one-to-one Boolean morphism.
\subsection{Proof of $\bullet_f$}\label{Proof:idempot}
By definition, the result holds true for any $A\in M_n\setminus\{\emptyset,\Omega_n,b_n,\sim b_n\}$\,.
It is also true for $A=\emptyset$ or $A=\Omega_n$\,, since $f_{n+1}\bigl(\mu_n(B),\mu_n(\emptyset)\bigl)=f_{n+1}\bigl(\mu_n(B),\emptyset\bigl)=\mu_n(B)=\mu_n\bigl(f_n(B,\emptyset)\bigr)$
and similarly
$f_{n+1}\bigl(\mu_n(B),\mu_n(\Omega_n)\bigl)=f_{n+1}\bigl(\mu_n(B),\Omega_{n+1}\bigl)=\mu_n(B)=\mu_n\bigl(f_n(B,\Omega_n)\bigr)$\,.
\\[5pt]
The true difficulties come from the cases $A=b_n$ or $A=\sim b_n$\,.\\
Subsequently, it is assumed $A=b_n$\,; the case $A=\sim b_n$ is quite similar.\\
It comes:
$$
f_{n+1}\bigl(\mu_n(B),\mu_n(b_n)\bigr)=(\mathrm{id}\cup T)
\biggl(\mu_n(B)\cap\Bigl(\bigcup_{i\in I_n}\bigl(\Pi_n(i)\times\Gamma_n(i)\bigr)\Bigr)\biggr)
=(\mathrm{id}\cup T)\Bigl(\bigcup_{i\in I_n}\bigl(B\cap\Pi_n(i)\bigr)\times\Gamma_n(i)\Bigr)\;.
$$
The existence of $f_n(B,b_n)$ necessary implies the case 0\,,
and there is $C\in M_{\nu+1}$ such that $B=C_{[n]}$\,.\\
By recursion hypotheses $\bullet_\mu$, it comes $B\cap\omega_{[n]}=(C\cap\omega)_{[n]}=\omega_{[n]}$ if $\omega\in C$\,, $=\emptyset$ if $\omega\in C$\,.\\
Moreover, $f_n(\omega_{[n]},b_n)\cap f_n(B,b_n)=f_n(\omega_{[n]}\cap B,b_n)=f_n(\omega_{[n]},b_n)$ if $\omega\in C$\,, $=\emptyset$ if $\omega\in C$\,.\\
As a consequence $\bigcup_{i\in I_n}\bigl(B\cap\Pi_n(i)\bigr)\times\Gamma_n(i)=\bigcup_{i\in I_n}\Pi_n(i)\times\bigl(f_n(B,b_n)\cap\Gamma_n(i)\bigr)$\,.
\\[3pt]
By $\tilde\beta3$,
$B\cap\omega_{[n]}=B\cap b_n\cap\omega_{[n]}=f_n(B,b_n)\cap b_n\cap\omega_{[n]}=f_n(B,b_n)\cap\omega_{[n]}$ for any $\omega\in\mu_\nu(b_\nu)$\,.\\
As a consequence $\bigcup_{i\in I_n}\bigl(B\cap\Pi_n(i)\bigr)\times\Gamma_n(i)=\bigcup_{i\in I_n}\bigl(f_n(B,b_n)\cap\Pi_n(i)\bigr)\times\Gamma_n(i)$\,.
\\
By applying the both results, it comes:
$$
f_{n+1}\bigl(\mu_n(B),\mu_n(b_n)\bigr)=
\Bigl(\bigcup_{i\in I_n}\bigl(f_n(B,b_n)\cap\Pi_n(i)\bigr)\times\Gamma_n(i)\Bigr)
\cup
\Bigl(\bigcup_{i\in I_n}\bigl(f_n(B,b_n)\cap\Gamma_n(i)\bigr)\times\Pi_n(i)\Bigr)\;.
$$
And by definition, $f_{n+1}\bigl(\mu_n(B),\mu_n(b_n)\bigr)=\mu_n\bigl(f_{n+1}(B,b_n)\bigr)$\,.
\subsection{Proof of $\tilde\beta1$}
For $A\not\in \{\mu_n(b_n),\sim \mu_n(b_n),\emptyset,\Omega_{n+1}\}$, the propriety is inherited from $n$ by applying $\bullet_f$\,.\\
The property is also obvious for $A\in \{\emptyset,\Omega_{n+1}\}$\,.
\\
The difficulty comes from $A=\mu_n(b_n)$ or $A=\sim \mu_n(b_n)$\,;
then notice that $A\ne\emptyset$ by construction.
\\[5pt]
It is now hypothesized $A=\mu_n(b_n)$ and $\mu_n(b_n)\subset B$\,;
the case $A=\sim \mu_n(b_n)$ is quite similar.\\
Then $T\bigl(\mu_n(b_n)\bigr)\subset T(B)$ and by lemma, corollary 1\&2:
$$
f_{n+1}\bigl(B,\mu_n(b_n)\bigr)=
\bigl(B\cap\mu_n(b_n)\bigr)\cup\bigl(T(B)\cap\sim\mu_n(b_n)\bigr)=\mu_n(b_n)\cup\sim\mu_n(b_n)
=\Omega_{n+1}\;.
$$
\subsection{Proof of $\tilde\beta2$}
For $A\not\in \{\mu_n(b_n),\sim \mu_n(b_n),\emptyset,\Omega_{n+1}\}$, the propriety is inherited from $n$ by applying $\bullet_f$\,.\\
The property is also obvious for $A\in \{\emptyset,\Omega_{n+1}\}$\,.
\\
The property is then immediate for $A\in\{\mu_n(b_n),\sim \mu_n(b_n)\}$\,, since $T(B_1\cup B_2)=T(B_1)\cup T(B_2)$.
\subsection{Proof of $\tilde\beta3$}
For $A\not\in \{\mu_n(b_n),\sim \mu_n(b_n),\emptyset,\Omega_{n+1}\}$, the propriety is inherited from $n$ by applying $\bullet_f$\,.
The property is also obvious for $A\in \{\emptyset,\Omega_{n+1}\}$\,.
\\
The difficulty comes from $A\in \{\mu_n(b_n),\sim \mu_n(b_n)\}$.
\\[5pt]
It is now hypothesized $A=\mu_n(b_n)$\,;
the case $A=\sim \mu_n(b_n)$ is quite similar.\\
The result is immediate from corollary 2 of lemma.
\subsection{Proof of $\tilde\beta4$}
For $A\not\in \{\mu_n(b_n),\sim \mu_n(b_n),\emptyset,\Omega_{n+1}\}$, the propriety is inherited from $n$ by applying $\bullet_f$\,.
The property is also obvious for $A\in \{\emptyset,\Omega_{n+1}\}$\,.
\\
The difficulty comes from $A\in \{\mu_n(b_n),\sim \mu_n(b_n)\}$\,.\\[5pt]
It is now hypothesized $A=\mu_n(b_n)$\,;
the case $A=\sim \mu_n(b_n)$ is quite similar.\\
By corollary 2 of lemma, $f_{n+1}\bigl(\sim B,\mu_{n}(b_n)\bigr)=\bigl(\mu_n(b_n)\setminus B\bigr)\cup\bigl(\sim\mu_n(b_n)\setminus T(B)\bigr)$\,.\\
By $\ell6$, $f_{n+1}\bigl(\sim B,\mu_{n}(b_n)\bigr)=
\sim\Bigl(\bigl(B\cap\mu_n(b_n)\bigr)\cup\bigl(T(B)\cap\sim\mu_n(b_n)\bigr)\Bigr)=\sim f_{n+1}\bigl(B,\mu_{n}(b_n)\bigr)$\,.
\subsection{Lemma 2.}
Let $C\in M_{n+1}$\,. Then:
$$
f_{n+1}\Bigl(f_{n+1}\bigl(C,\mu_n(b_n)\bigr),\mu_n(b_n)\Bigr)=
f_{n+1}\Bigl(f_{n+1}\bigl(C,\mu_n(b_n)\bigr),\sim\mu_n(b_n)\Bigr)=
f_{n+1}\bigl(C,\mu_n(b_n)\bigr)
$$
and
$$
f_{n+1}\Bigl(f_{n+1}\bigl(C,\sim\mu_n(b_n)\bigr),\mu_n(b_n)\Bigr)=
f_{n+1}\Bigl(f_{n+1}\bigl(C,\sim\mu_n(b_n)\bigr),\sim\mu_n(b_n)\Bigr)=
f_{n+1}\bigl(C,\sim\mu_n(b_n)\bigr)
\;.
$$
\begin{description}
\item[Proof.]
The result is derived for $f_{n+1}\bigl(C,\mu_n(b_n)\bigr)$\,; it is quite similar for $f_{n+1}\bigl(C,\sim\mu_n(b_n)\bigr)$\,.
\\[3pt]
By corollary 2 of lemma, $f_{n+1}\bigl(C,\mu_n(b_n)\bigr)=\bigl(C\cap\mu_n(b_n)\bigr)\cup\bigl(T(C)\cap\sim\mu_n(b_n)\bigr)$\,.\\
Since $T\bigl(f_{n+1}\bigl(C,\mu_n(b_n)\bigr)\bigr)=f_{n+1}\bigl(C,\mu_n(b_n)\bigr)$ by definition, the proof is done by corollary 2.
\item[$\Box\Box\Box$]\rien
\end{description}
\emph{Corollary.}
As a direct consequence, $f_{n+1}\bigl(f_{n+1}(B,A),A\bigr)=f_{n+1}\bigl(f_{n+1}(B,A),\sim A\bigr)=f_{n+1}(B,A)$\,, whenever $f_{n+1}(B,A)$ exists.
\subsection{Proof of $\tilde\beta6w$}
Assume $f_{n+1}(B,A)$ and $f_{n+1}(B,\sim A)$ exist and $f_{n+1}(B,A)=B$\,.\\
By lemma 2, $f_{n+1}(B,\sim A)=f_{n+1}\bigl(f_{n+1}(B,A),\sim A\big)=f_{n+1}(B,A)=B$\,.
\section{Proof: partial completeness}
\label{Appendix:ProofOfAlmostCompletude}
To be proved:\\[5pt]
Let $\phi\in\mathcal{L}$ be constructed without $\Box$ or $\Diamond$\,.
Then $\vdash\phi$ in DmBL$_\ast$ if and only if $H_{\mathcal{B}}(\phi)=\Omega_\infty$\,.
\\[5pt]
From now on, let $\mathcal{L}_b=\bigl\{\phi\in\mathcal{L}\,/\, \phi\mbox{ is constructed without }\Box\mbox{ or }\Diamond\bigr\}$\,.\\
In fact, it will be proved:
\begin{equation}
\label{the:Property}
\mbox{$(H_{\mathcal{B}})_\equiv$ is a Boolean isomorphism between $(\mathcal{L}_b)_\equiv$ and $M_\infty$\,,}
\end{equation}
where $(\mathcal{L}_b)_\equiv$ is the set of equivalence classes of $\mathcal{L}_b$ and $(H_{\mathcal{B}})_\equiv$ is inferred from $H_{\mathcal{B}}$\,.
\\[5pt]
The proof is based on a recursive construction of $\mathcal{L}_b$ similar to the definition of $M_\infty$\,.
\paragraph{Construction.}
Assume the sequence $(\Omega_n,M_n,h_n,f_n,\Lambda_n)_{n\in\Nset}$ being constructed.\\
The sequence $(L_n)_{n\in\Nset}$ is defined by:
\begin{itemize}
\item $L_0=\mathcal{L}_C$\,,
\item $L_{n+1}\subset\mathcal{L}_b$ is the set generated by $L_n$, the classical operators, the conditionals $(\cdot|\phi)$ where $\phi\in L_n$ and $H_{\mathcal{B}}(\phi)=b_n$\,.
\end{itemize}
The set $\Sigma_n\subset (L_n)_\equiv$ is defined as the generating partition of $L_n$,
\\
that is $\forall\phi\in (L_n)_\equiv\,,\;\exists S\subset\Sigma_n\,,\;\bigvee_{\sigma\in S}\sigma=\phi$ and $\sigma\wedge\sigma'=\bot$ for any $\sigma,\sigma'\in\Sigma_n$ such that $\sigma\ne\sigma'$\,.
\\[5pt]
The following property is proved in next paragraph:
\begin{equation}
\label{the:Equation}
\mathrm{card}(\Sigma_n)\le\mathrm{card}(\Omega_n)\;.
\end{equation}
Since $(H_{\mathcal{B}})_\equiv$ is by construction an onto morphism from $(L_n)_\equiv$ to $M_{n:\infty}$\,,
(\ref{the:Equation}) implies that $(H_{\mathcal{B}})_\equiv$ is a Boolean isomorphism between $(L_n)_\equiv$ and $M_{n:\infty}$\,.\\
These isomorphisms and the cyclic definition of $\Lambda_n$ then imply $\mathcal{L}_b=\cup_{n\in\Nset}(L_n)$\,.\\
As a conclusion, (\ref{the:Property}) is deduced.
\paragraph{Proof of (\ref{the:Equation}) for $n=0$.}
It is obvious, since $(M_{0:\infty},H_{\mathcal{B}})$ is a complete model for $\mathcal{L}_C$\,.
\paragraph{True for $n$ implies true for $n+1$.}\rien\\
The recursion hypothesis implies that $(H_{\mathcal{B}})_\equiv$ is an isomorphism between $(L_n)_\equiv$ and $M_{n:\infty}$\,.\\
Define then $\beta_n\in (L_n)_\equiv$ such that $(H_{\mathcal{B}})_\equiv(\beta_n)=b_n$\,.\\
It is known that $\bigl((\cdot|\beta_n)\big|\neg\beta_n\bigr)=(\cdot|\beta_n)$ and
$\bigl((\cdot|\neg\beta_n)\big|\beta_n\bigr)=(\cdot|\neg\beta_n)$\,.\\
Then, since sub-universes are classical, $\Sigma_{n+1}=\bigl\{\sigma\wedge(\sigma'|\beta_n)\wedge(\sigma''|\neg\beta_n)\;/\;\sigma,\sigma',\sigma''\in\Sigma_n\bigr\}\setminus\{\bot\}$\,.\\[5pt]
Now, denote $B_n
=\bigl\{\sigma\in\Sigma_n\,\big/\,\sigma\wedge \beta_n=\sigma\bigr\}$
and
$\overline{B}_n
=\bigl\{\sigma\in\Sigma_n\,\big/\,\sigma\wedge\neg\beta_n=\sigma\bigr\}$\,.\\
It comes that $(\sigma'|\beta_n)=(\sigma''|\neg\beta_n)=\bot$ for $\sigma'\in\overline{B}_n$ and $\sigma''\in B_n$\,.\\
Moreover $\sigma\wedge(\sigma'|\beta_n)\wedge(\sigma''|\neg\beta_n)=\bot$ for $\sigma\not\in \{\sigma',\sigma''\}$\,;
in particular
$\sigma\wedge(\sigma|\beta_n)=\sigma$ for $\sigma\in B_n$\,, and
$\sigma\wedge(\sigma|\neg\beta_n)=\sigma$ for $\sigma\in \overline{B}_n$\,.
\\[5pt]
Then, the two cases are considered:
\subparagraph{Case 1.}
$\Sigma_{n+1}=
\bigcup_{\sigma\in B_n}
\bigcup_{\sigma'\in \overline{B}_n}\bigl\{\sigma\wedge(\sigma'|\neg\beta_n),\sigma'\wedge(\sigma|\beta_n)\bigr\}
$\,, owing to above discussion;
then: $$
\mathrm{card}(\Sigma_{n+1})\le2\mathrm{card}(B_n)\mathrm{card}(\overline{B}_n)=2\mathrm{card}(b_n)\mathrm{card}(\sim b_n)=\mathrm{card}(\Omega_{n+1})\;.
$$
\subparagraph{Case 0.}
In this case, $\beta_n=\beta_\nu$\,.\\
Define $C_\nu=\{\sigma\in\Sigma_{\nu+1}/\sigma\wedge\beta_\nu=\sigma\}$
and $\overline{C}_\nu=\{\sigma\in\Sigma_{\nu+1}/\sigma\wedge\neg\beta_\nu=\sigma\}$\,.\\
Define also $D[\phi]=\{\sigma\in\Sigma_n/\sigma\wedge\phi=\sigma\}$ for any $\phi\in(L_{\nu+1})_\equiv$\,.\\[5pt]
From previously, it is know that $\Sigma_{n+1}$ contains elements of the form $\sigma\wedge(\sigma'|\neg\beta_n)$ or $\sigma'\wedge(\sigma|\beta_n)$ with $(\sigma,\sigma')\in B_n\times \overline{B}_n$\,;
but construction at step $\nu+1$ implies additional constraints, to be specified.
\\[3pt]
Let consider especially the case $\sigma\wedge(\sigma'|\neg\beta_n)$; case $\sigma'\wedge(\sigma|\beta_n)$ is quite similar.\\
Noticed that there is $(\tau_\sigma,\tau_\sigma')\in C_\nu\times\overline{C}_\nu$
such that $\sigma\in D[\tau_\sigma]\cap D[(\tau_\sigma'|\neg\beta_\nu)]$\,,
and $(\theta_{\sigma'},\theta_{\sigma'}')\in \overline{C}_\nu\times C_\nu$
such that $\sigma\in D[\theta_{\sigma'}]\cap D[(\theta_{\sigma'}'|\beta_\nu)]$\,.\\
Now, $\tau_\sigma\wedge(\tau_\sigma'|\neg\beta_\nu)\wedge\bigl(\theta_{\sigma'}\wedge(\theta_{\sigma'}'|\beta_\nu)\big|\neg\beta_\nu\bigr)=
\bigl(\tau_\sigma\wedge(\theta_{\sigma'}'|\beta_\nu)\bigr)\wedge
(\tau_\sigma'\wedge\theta_{\sigma'}|\neg\beta_\nu)=\bot$
unless $\tau_\sigma=\theta_{\sigma'}'$ and $\tau_\sigma'=\theta_{\sigma'}$\,.\\[3pt]
Finally, $\exists(\tau,\theta)\in C_\nu\times\overline{C}_\nu\,,\;(\sigma,\sigma')\in \bigl(D[\tau]\cap D[(\theta|\neg\beta_\nu)]\bigr)\times\bigl(D[\theta]\cap D[(\tau|\beta_\nu)]\bigr)$\,;
similarly the case $\sigma'\wedge(\sigma|\beta_n)$ implies $\exists(\theta,\tau)\in \overline{C}_\nu\times C_\nu\,,\;(\sigma,\sigma')\in \bigl(D[\theta]\cap D[(\tau|\neg\beta_\nu)]\bigr)\times\bigl(D[\tau]\cap D[(\theta|\neg\beta_\nu)]\bigr)$\,.\\[5pt]
At last $\mathrm{card}(\Sigma_{n+1})\le\sum_{(\tau,\theta)\in C_\nu\times\overline{C}_\nu} 2\,\mathrm{card}\bigl(D[\tau]\cap D[(\theta|\neg\beta_\nu)]\bigr) \mathrm{card}\bigl(D[\theta]\cap D[(\tau|\beta_\nu)]\bigr)$\\
\rien\hspace{175pt}$=\sum_{i\in I_n}2\,\mathrm{card}\bigl(\Pi_n(i)\bigr) \mathrm{card}\bigl(\Gamma_n(i)\bigr)=\mathrm{card}(\Omega_{n+1})\;.$
\section{Probability construction}
\label{Appendix:Probabilition}
To be proved:\\[5pt]
Let $\pi$ be a probability defined over $C$\,, such that $\pi(\phi)>0$ for any $\phi\not\equiv_C\bot$.
Then, there is a (multiplicative) probability $P$ defined over DmBL$_\ast$ such that  $\forall\phi\in\mathcal{L}_C\,,\;P(\phi)=\pi(\phi)$\,.
\\[5pt]
The construction is a recursion based on the definition of $\mathcal{B}$\,.
\subsection{Construction}
The probabilities $P_n|_{n\in\Nset}$ are defined over $M_{n:\infty}$ by:
$$
P_n(A_\infty)=\sum_{\omega\in A}P_n\bigl(\{\omega\}_\infty\bigr)\quad\mbox{for any }A\in M_n\,,
$$
and:
\subparagraph{Initialization.}\rien\\
For $\omega=\bigl(\delta_\theta|_{\theta\in\Theta}\bigr)\in\Omega_0$ and $\tau\in\Theta$, define $\tau_\omega=\tau$ if $\delta_\tau=1$ and $\tau_\omega=\neg\tau$ if $\delta_\tau=0$\,.\\[3pt]
Then set $P_0\bigl(\omega_\infty\bigr)=\pi\left(\bigwedge_{\theta\in\Theta}\theta_\omega\right)$ for any $\omega\in\Omega_0$\,,
\subparagraph{From n to n+1.}
For any $(\omega,\omega')\in\Pi_n(i)\times\Gamma_n(i)$\,, set:
$$
P_{n+1}\bigl((\omega,\omega')_\infty\bigr)=\frac{P_n(\omega_\infty)P_n(\omega'_\infty)}{P_n\bigl(\Gamma_n(i)_\infty\bigr)}
\quad\mbox{and}\quad
P_{n+1}\bigl((\omega',\omega)_\infty\bigr)=\frac{P_n(\omega_\infty)P_n(\omega'_\infty)}{P_n\bigl(\Pi_n(i)_\infty\bigr)}\;.
$$
An example of construction is given in appendix~\ref{BayesMod:Comp&Ex}.
\subparagraph{Notation.}
For $m\le n$ and $A\in M_m$, the probability $P_n(A_\infty)$ is denoted $P_n(A)$ for simplicity.
\subsection{Properties.}
\subsubsection{Proposition 1} $P_n\subset P_{n+1}$\,, \emph{i.e.} $P_{n+1}(A)=P_{n}(A)$ for any $A\in M_n$\,.
\begin{description}
\item[Proof.]
For $A\in M_n$\,,
$\mu_n(A)=\bigcup_{i\in I_n}\biggl(\Bigl(\bigl(A\cap\Pi_n(i)\bigr)\times\Gamma_n(i)\Bigr)\cup
\Bigl(\bigl(A\cap\Gamma_n(i)\bigr)\times\Pi_n(i)\Bigr)\biggr)$ and then:
$$\rien\hspace{-20pt}\begin{array}{@{}l@{}}\displaystyle
P_{n+1}(A)=\sum_{i\in I_n}\left(
\sum_{\omega\in A\cap\Pi_n(i)}\left(\sum_{\omega'\in\Gamma_n(i)}
\frac{P_n(\omega)P_n(\omega')}{P_n\bigl(\Gamma_n(i)\bigr)}\right)
+
\sum_{\omega\in A\cap\Gamma_n(i)}\left(\sum_{\omega'\in\Pi_n(i)}
\frac{P_n(\omega)P_n(\omega')}{P_n\bigl(\Pi_n(i)\bigr)}\right)
\right)
\vspace{4pt}\\\displaystyle
\rien\hspace{100pt}
=\sum_{i\in I_n}\left(
\sum_{\omega\in A\cap\Pi_n(i)}P_n(\omega)+
\sum_{\omega\in A\cap\Gamma_n(i)}P_n(\omega)
\right)=
\sum_{\omega\in A}P_n(\omega)=P_n(A)\;.
\end{array}$$
\item[$\Box\Box\Box$]\rien
\end{description}
\emph{Corollary.}
$P_n(\Omega_\infty)=1$\,.
\\[3pt]
Derived from $P_0(\Omega_\infty)=\pi(\top)=1$ which is obvious.
\\[4pt]
\emph{Corollary of the corollary.}
$P_n$ is indeed a probability in the classical meaning.
\\[3pt]
Additivity, coherence are obtained by construction.
Finiteness comes from the corollary.
\subsubsection{Proposition 2}
\label{prop:3}
\begin{enumerate}
\item $\displaystyle P_n\bigl(\Pi_n(i)\bigr)+P_n\bigl(\Gamma_n(i)\bigr)=\frac{P_n\bigl(\Pi_n(i)\bigr)}{P_n(b_n)}=\frac{P_n\bigl(\Gamma_n(i)\bigr)}{P_n(\sim b_n)}$\,,
for any $i\in I_n$\,, 
\item
$P_{n+1}\bigl(\mu_n(b_n)\cap A\bigr)=P_{n+1}(b_n)P_{n+1}\Bigl(f_{n+1}\bigl(A,\mu_n(b_n)\bigr)\Bigr)$\,, for any $A\in M_{n+1}$\,,
\item
$P_{n+1}\bigl(\sim \mu_n(b_n)\cap A\bigr)=P_{n+1}(\sim b_n)P_{n+1}\Bigl(f_{n+1}\bigl(A,\sim \mu_n(b_n)\bigr)\Bigr)$\,, for any $A\in M_{n+1}$\,.
\end{enumerate}
These propositions are proved recursively.
\begin{description}
\item[Proof of 1.]
Obvious in case 1; the difficulty arises for case 0.
\\[3pt]
Assume now case 0, and let $(\omega,\omega')\in I_n$\,, \emph{i.e.} $\omega\in\mu_\nu(b_\nu)$ and $\omega'\in\sim\mu_\nu(b_\nu)$.\\
Then $\frac{P_n\Bigl(f_{\nu+1}\bigl(\omega',\sim \mu_\nu(b_\nu)\bigr)\cap\omega\Bigr)}{P_n(b_\nu)}=P_n\biggl(f_{\nu+1}\Bigl(f_{\nu+1}\bigl(\omega',\sim \mu_\nu(b_\nu)\bigr)\cap\omega,\mu_\nu(b_\nu)\Bigr)\biggr)$, by the recursion hypothesis over 2,
and finally $\frac{P_n\bigl(\Pi_n(i)\bigr)}{P_n(b_n)}=P_n\Bigl(f_{\nu+1}\bigl(\omega',\sim \mu_\nu(b_\nu)\bigr)\cap f_{\nu+1}\bigl(\omega,\mu_\nu(b_\nu)\bigr)\Bigr)$\,.\\
Similarly, it is derived $\frac{P_n\bigl(\Gamma_n(i)\bigr)}{P_n(\sim b_n)}=P_n\Bigl(f_{\nu+1}\bigl(\omega',\sim \mu_\nu(b_\nu)\bigr)\cap f_{\nu+1}\bigl(\omega,\mu_\nu(b_\nu)\bigr)\Bigr)$\,.\\
Then $\frac{P_n\bigl(\Pi_n(i)\bigr)}{P_n(b_n)}=\frac{P_n\bigl(\Gamma_n(i)\bigr)}{P_n(\sim b_n)}$ and the result is deduced from $P_n(b_n)+P_n(\sim b_n)=1$.
\item[Proof of 2.]
Since $P_{n+1}\bigl(T(\omega)\bigr)=P_{n+1}(\omega)\frac{P_n\bigl(\Gamma_n(i)\bigr)}{P_n\bigl(\Pi_n(i)\bigr)}$ for $\displaystyle\omega\in \bigcup_{i\in I_n}\bigl(\Pi_n(i)\times\Gamma_n(i)\bigr)$\,,
it comes:
$$\begin{array}{@{}l@{}}
\displaystyle P_{n+1}\Bigl(f_{n+1}\bigl(A,\mu_n(b_n)\bigr)\Bigr)=\sum_{i\in I_n}\biggl(\sum_{\omega\in A\cap(\Pi_n(i)\times\Gamma_n(i))}P_{n+1}(\omega)\frac{P_n\bigl(\Pi_n(i)\bigr)+P_n\bigl(\Gamma_n(i)\bigr)}{P_n\bigl(\Pi_n(i)\bigr)}\biggr)
\vspace{4pt}\\\displaystyle
\rien\hspace{50pt}=\sum_{i\in I_n}\frac{P_{n+1}\Bigl(A\cap\bigl(\Pi_n(i)\times\Gamma_n(i)\bigr)\Bigr)}{P_n(b_n)}
=\frac{P_{n+1}\bigl(\mu_n(b_n)\cap A\bigr)}{P_n(b_n)}\;.
\end{array}$$
\item[Proof of 3.] Similar to 2.
\item[$\Box\Box\Box$]\rien
\end{description}
\subsubsection{Conclusion.}
\label{AppC:Conclude}
Define $P_\infty=\bigcup_{n\in\Nset}P_n$\,, that is $\forall A\in M_{n:\infty}\,,\;P_\infty(A)=P_n(A)$\,.
\\
By inheritance from $P_n$, $P_\infty$ is a probability over $M_{\infty}$, which verifies the property:
$$
P_{\infty}(A\cap B)=P_{\infty}(A)P_{\infty}\bigl(f_{\infty}(B,A)\bigr)\,,\quad\mbox{for any }A\in M_{\infty}\;.
$$
Define $P(\phi)=P_\infty\bigl(H_{\mathcal{B}}(\phi)\bigr)$\,\\
Since $\mathcal{B}$ is a model of DmBL, it is clear that the additivity, coherence and finiteness are inherited from $P_\infty$\,.\\[5pt]
Now, it is noticed that $\vdash\phi\times\psi$ is equivalent to $(\psi|\phi)\equiv\psi$\,.\\
It is then implied:
$$
P(\phi\wedge\psi)=P_\infty\bigl(H_{\mathcal{B}}(\phi)\cap H_{\mathcal{B}}(\psi)\bigr)= P_\infty\bigl(H_{\mathcal{B}}(\phi)\bigr)P_\infty\Bigl(f_{\infty}\bigl(H_{\mathcal{B}}(\psi),H_{\mathcal{B}}(\phi)\bigr)\Bigr)=
P_\infty\bigl(H_{\mathcal{B}}(\phi)\bigr)P_\infty\bigl(H_{\mathcal{B}}(\psi)\bigr)\,.
$$
Finally, $P$ verifies the multiplicativity.
\\[5pt]
\underline{$P$ is a probability over DmBL$_\ast$\,.}
\paragraph{Rational structure of $P$\,.}
Let $\Sigma=\left\{\left.\bigwedge_{\theta\in\Theta}\epsilon_\theta\;\right/\;\epsilon\in\prod_{\theta\in\Theta}\{\theta,\neg\theta\}\right\}$\,.
For any $\phi\in\mathcal{L}$\,, there is a rational function $R_\phi: \Rset^{(2^\Theta)}\rightarrow\Rset$ such that $P(\phi)=R_\phi\bigl(\pi(\sigma)|_{\sigma\in\Sigma}\bigr)$\,.
\\[5pt]
The proof is obvious from the construction.
\section{Conditional model: construction and comparison}
\label{BayesMod:Comp&Ex}
\paragraph{First step construction (examples).}
In this paragraph, the objects $\Omega_{k}, f_{k},\Lambda_{k},\mu_{k-1}|_{k=0,1}$\,, \emph{i.e.} one iteration,  are explicitly constructed, as well as the associated probability extensions $P_0, P_1$ (\emph{c.f.} appendix~\ref{Appendix:Probabilition}).
It is assumed that $\Omega_0=\{a,b,c\}$\,.
This hypothesis cannot hold actually, but the case is sufficiently simple to be handled, and sufficiently complex to be illustrative.
By the way, only the case 1 of the construction is considered.
Case 0 is intractable in a true example.
\\[5pt]
In the sequel, the list $\Lambda_n$ is represented by the sequence $\lambda_n(k)|_{k=s_n}^{f_n-1}$\,.
For simplicity, $\mu_k(A)$ and $A$ will be identified.
\subparagraph{k=0.}
By definition, the list $\Lambda_0$ contains the elements of $\mathcal{P}(\Omega_0)\setminus\{\emptyset,\Omega_0\}$\,.
In this example, it is chosen $\Lambda_0=\{a,b\},c,\{b,c\},a,\{c,a\},b$ and $P_0(a)=0.2$\,, $P_0(b)=0.3$\,, $P_0(c)=0.5$\,. 
\subparagraph{k=1.}
Case 1 holds with $b_0=\{a,b\}$\,.
It comes $\Omega_1=\bigl\{(a,c),(b,c),(c,a),(c,b)\bigr\}$\,,
\\
$\mu_0(a)=(a,c)$\,, $\mu_0(b)=(b,c)$\,, $\mu_0(c)=\{(c,a),(c,b)\}$\,, 
\\[3pt]
$P_1(a)=P_1((a,c))=\frac{P_0(a)P_0(c)}{P_0(\Gamma_0(b_0))}=\frac{0.2\times0.5}{0.5}=0.2$\,, $P_1(b)=P_1((b,c))=\frac{0.3\times0.5}{0.5}=0.3$\,,
\\[3pt]
$P_1((c,a))=\frac{P_0(a)P_0(c)}{P_0(\Pi_0(b_0))}=\frac{0.2\times0.5}{0.2+0.3}=0.2$\,, $P_1((c,b))=\frac{0.3\times0.5}{0.2+0.3}=0.3$\,,
\\[3pt]
$f_1\bigl(a,\{a,b\}\bigr)=(\mathrm{id}\cup T)\bigl((a,c)\bigr)=\bigl\{(a,c),(c,a)\bigr\}$\,,
$f_1\bigl(b,\{a,b\}\bigr)=(\mathrm{id}\cup T)\bigl((b,c)\bigr)=\bigl\{(b,c),(c,b)\bigr\}$\,,\\
$f_1\bigl((c,a),c\bigr)=(\mathrm{id}\cup T)\bigl((c,a)\bigr)=\bigl\{(a,c),(c,a)\bigr\}$\,,
$f_1\bigl((c,b),c\bigr)=(\mathrm{id}\cup T)\bigl((c,b)\bigr)=\bigl\{(b,c),(c,b)\bigr\}$\,,
$f_1(c,c)=(\mathrm{id}\cup T)\bigl(\{(c,a),(c,b)\}\bigr)=\Omega_1$\,,
(other cases are obvious)
\\[3pt]
and $\Lambda_1=\{b,c\},a,\{c,a\},b\ +\ \mathcal{P}(\Omega_1)\setminus\mathcal{P}(\Omega_0)\ +\
\{a,b\},c$\,.
\paragraph{Comparison with the Product Space CEA.}
The conditional of PSCEA \cite{goodman2} is defined by:
$$
(S|X)=\bigl((S\cap X)\times U\times\dots\bigr)
\cup \bigl(X^c\times (S\cap X)\times U\times\dots\bigr)
\cup \bigl(X^c\times X^c\times (S\cap X)\times U\times\dots\bigr)\cup\dots\;,
$$
where $X$ and $S$ are unconditioned propositions and $U$ is the universe.\\
As a consequence, $(S|X)=\bigl((S\cap X)\times U\times\dots\bigr)\cup \bigl(X^c\times (S|X)\bigr)$\,.\\
By identifying $X$ with $X\times U\times U\dots$, as suggested by \cite{goodman2}, the PSCEA appears as a fix point for the equation $(S|X)=(S\cap X)\cup \bigl(X^c\times (S|X)\bigr)$\,.
\\[3pt]
Let compare this result with the case 1 of the model construction of DmBL$_\ast$
(the case 0 is just a constrained version of case 1, accordingly to the logical rules of DmBL$_\ast$).\\
Assuming $A\in M_n$, it comes $f_{n+1}\bigl(\mu_n(A),\mu_n(b_n)\bigr)=\bigl(\mu_n(A)\cap\mu_n(b_n)\bigr)\cup\bigl(\sim b_n\times (A\cap b_n)\bigr)$ and $\sim\mu_n(b_n)\cap f_{n+1}\bigl(\mu_n(A),\mu_n(b_n)\bigr)=\sim b_n\times (A\cap b_n)$\,.\\
Finally, it is deduced $f_{n+1}\bigl(\mu_n(A),\mu_n(b_n)\bigr)=\bigl(\mu_n(A)\cap\mu_n(b_n)\bigr)\cup\Bigl(\sim\mu_n(b_n)\cap f_{n+1}\bigl(\mu_n(A),\mu_n(b_n)\bigr)\Bigr)$\,.\\
By the way, it is recognized in both cases the logical fix point $(\psi|\phi)=(\phi\wedge\psi)\vee\bigl(\neg\phi\cap(\psi|\phi)\bigr)$\,.
\\[3pt]
Nevertheless, the similitude is limited.
In the case of PSCEA, the fix point is defined globally for any $X$ simultaneously.
In the case of DmBL$_\ast$, the iterations are individual (\emph{i.e.} related to $b_n$ and $\sim b_n$);
moreover, the conditional of DmBL$_\ast$ is defined for any propositions.
\\[5pt]
Then although a true comparison of PSCEA with DmBL needs a dedicated work, it seems possible to state that PSCEA is somewhat more constrained than DmBL$_\ast$ (\emph{i.e.} DmBL$_\ast$ is likely to specify the unconditioned level less than PSCEA).
In particular, the logical independence, specific to DmBL$_\ast$, will not infer significant consequences over the unconditioned level, as it is shown subsequently.
\subparagraph{Logical independence and unconditional propositions.}
A thorough examination of the definition $(S|X)$ shows that $(S|X)=S$ is equivalent to 
$S=\emptyset$ or $S=U$ or $X=U$\,.
\\[3pt]
Now, let consider $A$ and $B$ two unconditioned propositions in the model of DmBL$_\ast$.
The model construction implies that $(B|A)$ is an unconditioned proposition if and only if $A=\emptyset$ or $A=\Omega_\infty$ or $B=\emptyset$ or $A\subset B$ (\emph{c.f.} also Lewis result discussed in section~\ref{Section:Proba:DmBL})\,.
Then $(B|A)=B$ if and only if $A=\emptyset$ or $A=\Omega_\infty$ or $B=\emptyset$ or $B=\Omega_\infty$\,.
\\[3pt]
It is seen that there is a slight difference, since DmBL$_\ast$ considers the case of empty conditioning: whatever, it has been seen that the behavior of DmBL$_\ast$ was quite specific regarding to $(\cdot|\bot)$.
But essentially, both theories reduces the logical independence between unconditioned propositions to trivial cases.
\\[3pt]
As a conclusion, the notion of logical independence does not change the status of the unconditioned propositions: there is not a significant difference between DmBL$_\ast$ and PSCEA at the unconditioned level about this notion.
\end{document}